\documentclass[12pt]{jloganal}
 
 \usepackage[english]{babel}
 \usepackage{amssymb}
 \usepackage{amsmath}
 \usepackage{graphicx}
 \usepackage{bbm}
 \usepackage[right]{eurosym}
 \usepackage[T1]{fontenc}
 
 \def\vs{\vspace}

 \def\id{\mathrm{id}}
 
 \def\exp{\mathrm{exp}}

 \newtheoremstyle{newline}% name
 {\baselineskip}% Space above
 {\baselineskip}% Space below
 {\itshape}% Body font
 {}% Indent amount
 {\bfseries}% Theorem head font
 {}% Punctuation after theorem head
 {\newline}% Space after theorem head
 {\thmname{#1}\thmnumber{ #2}\thmnote{ (#3)}}% Theorem head spec (can be left empty, meaning ‘normal’ )
 
 \theoremstyle{newline}
 \newtheorem{thm}{Theorem}[section]
 \newtheorem{defn}[thm]{Definition}
 \newtheorem{rem}[thm]{Remark}
 \newtheorem{prop}[thm]{Proposition}
 \newtheorem{fac}[thm]{Fact}
 \newtheorem{exc}[thm]{Example}
 \newtheorem{cor}[thm]{Corollary}

\title{On Preparation Theorems for $\mathbb{R}_{\textnormal{an,exp}}$-definable Functions}

\author{Andre Opris}
\address{Andre Opris, University of Passau, 
	Faculty of Computer Science and Mathematics, D-94030 Germany}
\email{andre.opris@uni-passau.de}
\date{October, 28, 2021}


\begin{document}

\begin{abstract}
\textit{Abstract:} In this article we give strong versions for preparation theorems for $\mathbb{R}_{\textnormal{an,exp}}$-definable functions outgoing from methods of Lion and Rolin ($\mathbb{R}_{\textnormal{an,exp}}$ is the o-minimal structure generated by all restricted analytic functions and the global exponential function). By a deep model theoretic fact of Van den Dries, Macintyre and Marker every $\mathbb{R}_{\textnormal{an,exp}}$-definable function is piecewise given by $\mathcal{L}_{\textnormal{an}}(\exp,\log)$-terms where $\mathcal{L}_{\textnormal{an}}(\exp,\log)$ denotes the language of ordered rings augmented by all restricted analytic functions, the global exponential and the global logarithm. The idea is to consider log-analytic functions at first, i.e. functions which are iterated compositions from either side of globally subanalytic functions and the global logarithm, and then $\mathbb{R}_{\textnormal{an,exp}}$-definable functions as compositions of log-analytic functions and the global exponential.
\end{abstract}
\maketitle

\section*{Introduction}
This paper contributes to analysis in the framework of o-minimal structures. O-minimality is a concept from mathematical logic with connections and applications to geometry, analysis, number theory and other areas. Sets and functions definable in an o-minimal structure (i.e. ‘belonging to’) exhibit tame geometric and combinatorial behaviour. We refer to the book of Van den Dries~\cite{dries1998} for the general properties of o-minimal structures; in the preliminary section we state the definition and give examples.

Log-analytic functions have been defined by Lion and Rolin in their seminal paper \cite{Lion1997} (see also \cite{KaiserOpris} for the definition). They are iterated compositions from either side of globally subanalytic functions (see \cite{10.1215/S0012-7094-96-08416-1}) and the global logarithm. 

{\bf Example}\\
{\it The function 
$$f:\mathbb{R}^2 \to \mathbb{R}, (t,x) \mapsto \arctan(\log(\max\{\log(t^4+\log(x^2+2)),1\})),$$ 
is log-analytic.} 

Their definition is kind of hybrid as we will explain. The globally subanalytic functions
are precisely the functions that are definable in the o-minimal structure $\mathbb{R}_\textnormal{an}$ of restricted analytic functions (see \cite{10.1215/S0012-7094-96-08416-1}). Because the global logarithm is not globally subanalytic the class of log-analytic functions contains properly the class of $\mathbb{R}_{\textnormal{an}}$-definable functions. Since from the global logarithm the exponential function is definable and the latter is not log-analytic the class of log-analytic functions is not a class of definable functions. It is properly contained in the class of $\mathbb{R}_{\textnormal{an,exp}}$-definable functions ($\mathbb{R}_{\textnormal{an,exp}}$ is the o-minimal structure generated by all restricted analytic functions and the global exponential function, see \cite{10.1215/S0012-7094-96-08416-1}). Hence log-analytic functions capture $\mathbb{R}_\textnormal{an}$-definability but not full $\mathbb{R}_{\textnormal{an,exp}}$-definability. To obtain full $\mathbb{R}_{\textnormal{an,exp}}$-definability one has to consider compositions of log-analytic functions and the global exponential (see \cite{Dries1994TheET}).

Outgoing from this considerations Lion and Rolin found a way to prove that $\mathbb{R}_{\textnormal{an}}$ and $\mathbb{R}_{\textnormal{an,exp}}$ have quantifier elimination without using model theoretical techniques (see \cite{Lion1997}). The o-minimality of these
structures follow from this quantifier elimination result. Their main idea was to introduce and prove preparation theorems: once a definable function $f(t,x)$ is
cellwise prepared with respect to a given variable $x$, it is much easier to solve
the equation $f(t,x) = 0$ with respect to the unknown $x$. In some sense, these statements can be understood as a monomialization of the definable functions with respect to $x$,
in the spirit of resolution of singularities for functions which involve the logarithm and exponential functions. There are two mistakes in \cite{Lion1997}. The first one, which concerns the log-analytic statement, has been corrected in \cite{WiesPaw1999}.
The second one, which concerns the proof of the log-exp-analytic statement, has
been fixed in \cite{Dries2002OminimalPT}, in the framework of the exponential closure of polynomially bounded o-minimal structures (on the reals).\\
These preparation theorems are crucial to deal with questions from analysis in the o-minimal context: Van den Dries and Miller proved a parametric version of Tamm's theorem for globally subanalytic functions (see \cite{Tamm}) and 
Tobias Kaiser showed the impressive fact that a real analytic globally subanalytic function has a holomorphic extension which is again globally subanalytic (see \cite{KaiserComplex}). But it is not straightforward to generalize these results on $\mathbb{R}_{\text{an,exp}}$-definable functions: An easy counterexample shows that Tamm's Theorem does not hold in this setting and the second question is completely open in this framework. To obtain a positive answer for the second issue one idea would be to look closer at preparing $\mathbb{R}_{\textnormal{an,exp}}$-definable functions.

For the rest of the introduction ''definable'' means ''$\mathbb{R}_{\textnormal{an,exp}}$-definable'' .

%and the former is not true as easy counterexamples show. The idea is also to compute such a holomorphic extension To deal with such questions the idea is to look closer at preparing $\mathbb{R}_{\textnormal{an,exp}}$-definable functions. 
\iffalse
Furthermore it is hard to control the function on the boundary of the piece as the following example indicates.

{\bf Example}\\
{\it Let $f(t,x)$ be the efinable function with $f(t,x)=\exp(-1/x)-t/\log(x)$ if $x>0$ and zero otherwise. Then the following asymptotics hold. For every $t \neq 0$ we have $f(t,-) \sim -t / \log(x)$ and $f(0,-) \sim \exp(-1/x)$ as $x \searrow 0$.}
\fi
For example a log-analytic function $f(t,x)$ where $x$ is the last variable can be cellwise prepared as $f(t,x)=a(t)\vert{y_0(t,x)}\vert^{q_0} \cdot ... \cdot \vert{y_r(t,x)}\vert^{q_r}u(t,x)$ where $y_0(t,x)=x-\Theta_0(t)$ and inductively $y_l(t,x)=\log(\vert{y_{l-1}(t,x)}\vert)-\Theta_l(t)$ for $l \in \{1,...,r\}$, the $q_j$'s are rational exponents and $u(t,x)$ is a unit of a special form. This gives roughly that the function $f(t,-)$ behaves cellwise as iterated logarithms independently of $t$ where the order of iteration is bounded in terms of $f$ (compare also with Van den Dries and Speissegger in \cite{Dries2002OminimalPT}). Note that the functions $a(t),\Theta_0(t),...,\Theta_r(t)$ although being definable are in general not log-analytic anymore. We will present an example below. But it turns out that a log-analytic function can be prepared with data of a special form. \\
Given the projection $\pi:\mathbb{R}^n \times \mathbb{R} \to \mathbb{R}^n, (t,x) \mapsto t,$ on the first $n$ coordinates and a definable cell $C \subset \mathbb{R}^n \times \mathbb{R}$, we call a function $g:\pi(C) \to \mathbb{R}$ \textit{$C$-nice} if $g$ is the composition of log-analytic functions and exponentials of the form $\exp(h)$ where $h$ is the component of a center of a logarithmic scale on $C$ (compare with Definition 2.1 for the definition of a logarithmic scale, the logarithmic scale may depend on $g$). We call such a function $h$ \textit{$C$-heir} (see Definition \ref{2.29} below). It is immediately seen that a log-analytic function on $\pi(C)$ is $C$-nice and that every $C$-nice function is definable. Note that the class of $C$-nice functions does not necessarily coincide with the class of definable functions: if the cell $C$ is \textit{simple}, i.e. for every $t \in \pi(C)$ there is $d_t \in \mathbb{R}_{>0} \cup \{\infty\}$ such that $C_t= \textnormal{}]0,d_t[$ (see for example Definition 2.15 in \cite{KaiserOpris}), the class of $C$-nice functions coincides with the class of log-analytic ones (since in \cite{KaiserOpris} it is shown that the center of a logarithmic scale vanishes on such a cell). 
%eventuell: If $C$ is open one can show that a $C$-nice function is a composition of log-analytic functions and exponentials of locally bounded functions, i.e. a \textit{restricted log-exp-analytic function} as defined in \cite{WiesławPawłucki1999}. So one sees that the notion of $C$-nice functions is a much better description than only definability. 

With the class of $C$-nice functions we are able to prove the following preparation theorem for log-analytic functions. For practical reasons we consider finitely many log-analytic functions and prepare them simultaneously.

{\bf Theorem A }\\
{\it Let $X \subset \mathbb{R}^n \times \mathbb{R}$ be definable, $m \in \mathbb{N}$ and let $f_1,...,f_m:X \to \mathbb{R}$ be log-analytic. Then there is $r \in \mathbb{N}_0$ and a decomposition $\mathcal{C}$ of $X_{\neq 0}:=\{(t,x) \in X \mid x \neq 0\}$ into finitely many definable cells such that for $C \in \mathcal{C}$ there is $\Theta:=(\Theta_0,...,\Theta_r)$ such that the functions $f_1|_C,...,f_m|_C$ are nicely log-analytically prepared in $x$ with center $\Theta$, i.e. 
$f_j(t,x)=a_j(t)\vert{y_0(t,x)}\vert^{q_0} \cdot ... \cdot \vert{y_r(t,x)}\vert^{q_r}u_j(t,x)$ for $(t,x) \in C$ and $j \in \{1,...,m\}$ where $y_0(t,x)=x-\Theta_0(t)$ and inductively $y_l(t,x)=\log(\vert{y_{l-1}(t,x)}\vert)-\Theta_l(t)$ for $l \in \{1,...,r\}$, $q_0,...,q_r \in \mathbb{Q}$, the functions $a_j,\Theta_j$ on $\pi(C)$ are $C$-nice and $u_j$ is a unit of the form $v_j \circ \phi$ where $v_j$ is a real power series which converges on an open neighborhood of $[-1,1]^s$ and $\phi:=(\phi_1,...,\phi_s)$ with $\phi(C) \subset [-1,1]^s$ and $\phi_i(t)=b_i(t)\vert{y_0(t,x)}\vert^{p_{i0}} \cdot ... \cdot \vert{y_r(t,x)}\vert^{p_{ir}}$, where $b_i$ is $C$-nice and $p_{i0},...p_{ir} \in \mathbb{Q}$.} 

This is a more precise preparation theorem than the version from \cite{Lion1997} since the data of the preparation is $C$-nice (and not only definable).\\
For example if $C$ is simple we have that $a,\Theta_0,...,\Theta_r$ and $b_1,...,b_s$ are log-analytic and therefore the preparation keeps log-analyticity (which cannot be obtained by using the preparation theorem from \cite{Lion1997}). Consequently if $g(t):=\lim_{x \searrow 0} f(t,x)$ exists for a log-analytic $f: \mathbb{R}^n \times \mathbb{R} \to \mathbb{R}, (t,x) \mapsto f(t,x)$, we see that $g$ is also log-analytic (compare with Theorem 3.1 in \cite{KaiserOpris}). Outgoing from this consideration important differentiability properties for log-analytic functions like non-flatness or a parametric version of Tamm's theorem could be established in \cite{KaiserOpris}. (Note also that a special version of Theorem A has been proven there: This is Theorem 2.30, a preparation theorem for log-analytic functions on simple cells.)
%If $C$ is open a prepared log-analytic function is the composition of log-analytic functions and exponentials of locally bounded functions, i.e. is restricted log-exp-analytic.

Our second goal for this paper is to establish a preparation theorem for definable functions again by adapting the arguments of Lion-Rolin from \cite{Lion1997} by considering definable functions as compositions of log-analytic functions and the global exponential. For this we introduce the concept of the \textit{exponential number} of a definable function $f$ with respect to a set $E$ of positive definable functions (or exponentials), i.e. the minimal number of iterations of functions from $E$ which are necessary to write $f$ as a composition of log-analytic functions and functions from $E$. 

{\bf Example}\\
{\it The function $\mathbb{R} \to \mathbb{R}, x \mapsto \exp(x^2+1)$, has exponential number at most $1$ with respect to $E:=\{\exp(x^2+1)\}$. Since $\mathbb{R}_{>0} \to \mathbb{R}, y \mapsto \log(y+1),$ is log-analytic, the function 
$$f:\mathbb{R} \to \mathbb{R}, x \mapsto \log(\exp(x^2+1)+1),$$
has also exponential number at most 1 with respect to $E$. The function 
$$\mathbb{R} \to \mathbb{R}, x \mapsto \exp(\exp(x^2+1)),$$
has exponential number at most $2$ with respect to $\{\exp(x^2+1),\exp(\exp(x^2+1))\}$. 
}

Every function $f:X \to \mathbb{R}$ which has exponential number at most $e \in \mathbb{N}_0$ with respect to a set of positive definable functions on $X$ is definable and for every definable function $g:X \to \mathbb{R}$ there is $e \in \mathbb{N}_0$ and a set $E$ of positive definable functions on $X$ such that $g$ has exponential number at most $e$ with respect to $E$ (see Section 1). The notion of the exponential number allows us to prove the following preliminary version of the preparation theorem for definable functions which gives a connection between the functions contained in $E$ and the exponentials which occur in the cellwise preparation of $f$.

{\bf Theorem B}\\
{\it Let $e \in \mathbb{N}_0$ and $m \in \mathbb{N}$. Let $X \subset \mathbb{R}^m$ be definable and let $f:X \to \mathbb{R}$ be a function. Let $E$ be a set of positive definable functions on $X$ such that $f$ has exponential number at most $e \in \mathbb{N}_0$ with respect to $E$. Then there is a decomposition $\mathcal{D}$ of $X$ into finitely many definable cells such that for every $D \in \mathcal{D}$ we have that $f|_D=a \cdot \exp(c) \cdot u$ where $u$ is a unit of the form $v \circ \phi$ where $v$ is a real power series and $\phi:=(\phi_1,...,\phi_s)$ with $\phi_i(t)=b_i(t)\exp(d_i(t))$. Additionally $a$ resp. $b_i$ are log-analytic and $c$ resp. $d_i$ are finite $\mathbb{Q}$-linear combinations of functions from $\log(E):=\{\log(g) \mid g \in E\}$ which have exponential number at most $e-1$ with respect to $E$.}

Finally we combine Theorem A and Theorem B and prove Theorem C, a preparation theorem for definable functions which is the main result of this paper.

Let $C \subset \mathbb{R}^n \times \mathbb{R}_{\neq 0}$ be a definable cell, let $r \in \mathbb{N}_0$ and let $E$ be a set of positive definable functions on $C$. Let $\Theta_0$,...,$\Theta_r$ be $C$-nice. We call a function $f:C \to \mathbb{R}$ $(-1,r)$-prepared with center $\Theta:=(\Theta_0,...,\Theta_r)$ with respect to $E$ if $f$ is the zero function. For $e \in \mathbb{N}_0$ call a function $f:C \to \mathbb{R}$ $(e,r)$-prepared with center $\Theta$ with respect to $E$ if for $(t,x) \in C$
$$f(t,x) = a(t) \vert{y_0(t,x)}\vert^{q_0} \cdot ... \cdot \vert{y_r(t,x)}\vert^{q_r} \exp(c(t,x)) \cdot u(t,x)$$
where $q_0,...,q_r \in \mathbb{Q}$, $y_0=x-\Theta_0(t)$, inductively $y_i=\log(\vert{y_{i-1}}\vert)-\Theta_i(t)$ for $i \in \{1,...,r\}$, $a$ is $C$-nice, $\exp(c) \in E$, and the unit $u$ is of the form $u = v \circ \phi$ where $v$ is a real power series with $v(C) \subset [-1,1]^s$ and $\phi=(\phi_1,...,\phi_s)$ with $\phi_i=b_i(t)\vert{y_0(t,x)}\vert^{p_{i0}} \cdot ... \cdot \vert{y_r(t,x)}\vert^{p_{ir}} \cdot \exp(d_i(t,x))$ where, $\phi_i(C) \subset [-1,1]$, $b_i$ is $C$-nice, $p_{i0},...,p_{ir} \in \mathbb{Q}$ and $\exp(d_i) \in E$. Additionally $c$ and $d_i$ are $(e-1,r)$-prepared with center $\Theta$ with respect to $E$.

{\bf Theorem C }\\
{\it Let $e \in \mathbb{N}_0$. Let $X \subset \mathbb{R}^n \times \mathbb{R}$ be definable and let $E$ be a set of positive definable functions on $X$. Let $f:X \to \mathbb{R}$ be a function with exponential number at most $e$ with respect to $E$. Then there is $r \in \mathbb{N}_0$ and a definable cell decomposition $\mathcal{C}$ of $X_{\neq 0}$ such that for every $C \in \mathcal{C}$ there is a finite set $P$ of positive definable functions on $C$ and $\Theta:=(\Theta_0,...,\Theta_r)$ such that the function $f|_C$ is $(e,r)$-prepared with center $\Theta$ with respect to $P$ and the following holds.
%Einführung C-Erbe, C-nett informell erklären!
\begin{itemize}
	\item [(1)]
	For every $g \in \log(P)$ there is $m \in \{-1,...,e-1\}$ such that $g$ is $(m,r)$-prepared in $x$ with center $\Theta$ with respect to $P$.
	\item [(2)]
	If $g \in \log(P):=\{\log(h) \mid h \in P\}$ is $(l,r)$-prepared in $x$ with center $\Theta$ with respect to $P$ for $l \in \{-1,...,e-1\}$ then $g$ is a finite $\mathbb{Q}$-linear combination of functions from $\log(E)$ restricted to $C$ which have exponential number at most $l$ with respect to $E$.
\end{itemize}}

Theorem C gives a cellwise preparation of definable functions in terms of the exponential number and with a precise statement how the unit looks like (compare with \cite{Lion1997}, Section 2 for a version of Theorem C not in terms of the exponential number and without mentioning any unit $u$; compare with \cite{Dries2002OminimalPT}, Section 5 for a version of Theorem C in terms of the exponential number where a unit $u$ is mentioned, but without any precise statement how this unit looks like.) In contrast to the treatment in \cite{Dries2002OminimalPT} or \cite{Lion1997} there is also precise information about the exponentials which occur in the preparation of $f$ on every cell $C$: they are of the form $\exp(h)$ where $h$ is the component of a center of a logarithmic scale on $C$ (coming from Theorem A, the log-analytic preparation) or $h$ is itself prepared and a finite $\mathbb{Q}$-linear combination of functions from $\log(E)$ restricted to $C$ (coming from Theorem B). \\
All this improvements together allow us to investigate natural intermediate classes between log-analytic and definable functions and discuss various analytic properties of them. One example are the so-called \textit{restricted log-exp-analytic functions}. These are definable functions which are compositions of log-analytic functions and exponentials of locally bounded functions (see Definition 2.5 in \cite{10.1215/00192082-10234104}). The global exponential function for example is restricted log-exp-analytic. Outgoing from Theorem C the differentiability results for log-analytic functions from \cite{KaiserOpris} could be generalized to this class of functions (see \cite{10.1215/00192082-10234104}). Another consequence of Theorem C is the deep technical result that a real analytic restricted log-exp-analytic function has a holomorphic extension which is again restricted log-exp-analytic (see \cite{OprisComplex}).

This paper is organised as follows. After a preliminary section about o-minimality, notations and conventions we give in Section 1 the definition of log-analytic functions, the definition of the exponential number of a definable function with respect to a set of positive definable functions and give elementary properties. In Section 2 we present the preparation theorem of Lion-Rolin from $\cite{Lion1997}$ for log-analytic functions, introduce 'nice functions' and give a proof for Theorem A above. Section 3 is devoted to the proof of Theorem B and Theorem C.
\section*{Preliminaries}
\subsection*{O-Minimality}
We give the definition and examples of o-minimal structures.

{\bf Semialgebraic sets:}\\
A subset $A$ of $\mathbb{R}^n$, $n \geq 1$, is called \textit{semialgebraic} if there are $k,l \in \mathbb{N}_0$ and real polynomials $f_i,g_{i,1},...,g_{i,k} \in \mathbb{R}[X_1,...,X_n]$ for $1 \leq i \leq l$ such that 
$$A = \bigcup_{i=1}^l \{x \in \mathbb{R}^n \mid f_i(x)=0, g_{i,1}(x)>0,...,g_{i,k}(x)>0\}.$$
A map is called semialgebraic if its graph is semialgebraic.

{\bf Semi- and subanalytic sets:}\\
A subset $A$ of $\mathbb{R}^n$, $n \geq 1$, is called \textit{semianalytic} if for every $a \in \mathbb{R}^n$ there are open neighborhoods $U,V$ of $a$ with $\overline{U} \subset V$, $k,l \in \mathbb{N}_0$ and real analytic functions $f_i,g_{i,1},...,g_{i,k}$ on $V$ for $1 \leq i \leq l$, such that
$$A \cap U = \bigcup_{i=1}^l\{x \in U \mid f_i(x)=0, g_{i,1}(x)>0,...,g_{i,k}(x)>0\}.$$
A subset $B$ of $\mathbb{R}^n$, $n \geq 1$, is called \textit{subanalytic} if for every $a \in \mathbb{R}^n$ there is an open neighborhood $U$ of $a$, some $p \geq n$ and some bounded semianalytic set $A \subset \mathbb{R}^p$ such that $B \cap U=\pi_n(A)$ where $\pi_n:\mathbb{R}^p \to \mathbb{R}^n, (x_1,...,x_p) \mapsto (x_1,...,x_n),$ is the projection on the first $n$ coordinates.
A map is called semianalytic or subanalytic if its graph is a semianalytic resp. subanalytic set, respectively. A set is called \textit{globally semianalytic} or \textit{globally subanalytic} if it is semianalytic or subanalytic, respectively, in the ambient projective space (or equivalently, after applying the semialgebraic homeomorphism $\mathbb{R}^n \to \textnormal{}]-1,1[^n, x_i \mapsto x_i/\sqrt{1+x_i^2}$.)

{\bf O-minimal structures:}\\
A \textit{structure} on $\mathbb{R}$ is axiomatically defined as follows. For $n \in \mathbb{N}$ let $M_n$ be a set of subsets of $\mathbb{R}^n$ and let $\mathcal{M}:=(M_n)_{n \in \mathbb{N}}$. Then $\mathcal{M}$ is a structure on $\mathbb{R}$ if the following holds for all $m,n,p \in \mathbb{N}$.
\begin{itemize}
	\item [(S1)] If $A,B \in \mathbb{R}^n$ then $A \cup B$, $A \cap B$ and $\mathbb{R}^n \setminus A \in M_n$. (So $M_n$ is a Boolean algebra of subsets of $M_n$.)
	\item [(S2)] If $A \in M_n$ and $B \in M_m$ then $A \times B \in M_{n+m}$.
	\item [(S3)] If $A \in M_p$ and $p \geq n$ then $\pi_n(A) \in M_n$ where $\pi_n:\mathbb{R}^p \to \mathbb{R}^n, (x_1,...,x_p) \mapsto (x_1,...,x_n)$, denotes the projection on the first $n$ coordinates.
	\item [(S4)] $M_n$ contains the semialgebraic subsets of $\mathbb{R}^n$. 
\end{itemize}
The structure $\mathcal{M}=(M_n)_{n \in \mathbb{N}}$ on $\mathbb{R}$ is called \textit{o-minimal} if additionally the following holds.
\begin{itemize}
	\item [(O)] The sets in $M_1$ are exactly the finite unions of intervals and points.
\end{itemize}

A subset of $\mathbb{R}^n$ is called \textit{definable} in the structure $\mathcal{M}$ if it belongs to $M_n$.
A function is definable in $\mathcal{M}$ if its graph is definable in $\mathcal{M}$.
The o-minimality axiom (O) implies that a subset of $\mathbb{R}$ which is definable in an o-minimal structure on $\mathbb{R}$ has only finitely many connected components. But much more can be deduced from o-minimality. A definable subset of $\mathbb{R}^n$, $n \in \mathbb{N}$ arbitrary, has only finitely many connected components and these are again definable. More generally, sets and functions definable in an o-minimal
structure exhibit tame geometric behaviour, for example the existence of definable
cell decomposition. We refer to the book of Van den Dries \cite{dries1998} for this and more of the general properties of o-minimal structures.

{\bf Examples of o-minimal structures:}
\begin{itemize}
	\item [(1)] The smallest o-minimal structure on $\mathbb{R}$ is given by the semialgebraic sets. It is denoted by $\mathbb{R}$.
	\item [(2)] $\mathbb{R}_{\textnormal{exp}}$, the structure generated on the real field by the global exponential function $\exp:\mathbb{R} \to \mathbb{R}_{>0}$ (i.e. the smallest structure containing the semialgebraic sets and the graph of the global exponential function), is o-minimal.
	\item [(3)] $\mathbb{R}_{\textnormal{an}}$, the structure generated on the real field by the restricted analytic functions, is o-minimal. A function $f:\mathbb{R}^n \to \mathbb{R}$ is called restricted analytic if there is a function $g$ that is real analytic on a neighborhood of $[-1,1]^n$ such that $f=g$ on $[-1,1]^n$ and $f=0$ otherwise. The sets definable in $\mathbb{R}_{\textnormal{an}}$ are precisely the globally subanalytic ones.
	\item [(4)] $\mathbb{R}_{\textnormal{an,exp}}$, the structure generated by $\mathbb{R}_{\textnormal{an}}$ and $\mathbb{R}_{\textnormal{exp}}$, is o-minimal.
\end{itemize}

\subsection*{Notations}

The empty sum is by definition $0$ and the empty product is by definition $1$. By $\mathbb{N}=\{1,2,...\}$ we denote the set of natural numbers and by $\mathbb{N}_0=\{0,1,2,...\}$ the set of nonnegative integers. We set $\mathbb{R}_{>0}:=\{x \in \mathbb{R} \mid x > 0\}$. For $m,n \in \mathbb{N}$ we denote by $M(m \times n,\mathbb{Q})$ the set of $m \times n$-matrices with rational entries. Given $x \in \mathbb{R}$ let $\textnormal{sign}(x) \in \{\pm 1\}$ its sign. By $\log_k$ we denote the $k$-times iterated of the natural logarithm. By ''$\sim$'' we denote asymptotic equivalence.

For $m \in \mathbb{N}$ and a set $X \subset \mathbb{R}^m$ we set the following: For a set $E$ of positive real valued functions on $X$ we set $\log(E):=\{\log(g) \mid g \in E\}$.  For $C \subset \mathbb{R}^m$ with $C \subset X$ and a set $E$ of real valued functions on $X$ we set $E|_C:=\{g|_C \mid g \in E\}$. 

For $X \subset \mathbb{R}^n \times \mathbb{R}$ let $X_{\neq 0}:=\{(t,x) \in X \mid x \neq 0\}$ and for $t \in \mathbb{R}^n$ we set $X_t:=\{x \in \mathbb{R} \mid (t,x) \in X\}$.

\subsection*{Convention}
Definable means definable in $\mathbb{R}_{\textnormal{an,exp}}$ if not otherwise mentioned.

\section{Log-Analytic Functions and the Exponential Number}
We fix $m \in \mathbb{N}$ and a definable $X \subset \mathbb{R}^m$.

\begin{defn}\label{1.1} 
Let $f:X \to \mathbb{R}$ be a function.
\begin{itemize}
	\item [(a)] Let $r \in \mathbb{N}_0$. By induction on $r$ we define that $f$ is \textbf{log-analytic of order at most} $r$.\\
	
	\textbf{Base case}: The function $f$ is log-analytic of order at most $0$ if $f$ is piecewise the restriction of  globally subanalytic functions, i.e. there is a decomposition $\mathcal{C}$ of $X$ into finitely many definable cells such that for $C \in \mathcal{C}$ there is a globally subanalytic function $F:\mathbb{R}^m \to \mathbb{R}$ such that $f|_C = F|_C$.\\
	
	\textbf{Inductive step}: The function $f$ is log-analytic of order at most $r$ if the following holds: There is a decomposition $\mathcal{C}$ of $X$ into finitely many definable cells such that for $C \in \mathcal{C}$ there are $k,l \in \mathbb{N}_{0}$, a globally subanalytic function $F:\mathbb{R}^{k+l} \to \mathbb{R}$, and log-analytic functions $g_1,...,g_k:C \to \mathbb{R}, h_1,...,h_l:C \to \mathbb{R}_{>0}$ of order at most $r-1$ such that
	$$f|_C=F(g_1,...,g_k,\log(h_1),...,\log(h_l)).$$
	
	\item[(b)] Let $r \in \mathbb{N}_0$. We call $f$ \textbf{log-analytic of order} $r$ if $f$ is log-analytic of order at most $r$ but not of order at most $r-1$.
	
	\item[(c)] We call $f$ \textbf{log-analytic} if $f$ is log-analytic of order $r$ for some $r \in \mathbb{N}_0$.
\end{itemize}
\end{defn}

\begin{rem}\label{1.2}
The following properties hold.
\begin{itemize}
	\item[(1)] A log-analytic function is definable.
	\item[(2)] The log-analytic functions are precisely those definable functions which are piecewise given by $\mathcal{L}_{\textnormal{an}}(^{-1},(\sqrt[n]{...})_{n=2,3,...}, \log)$-terms (compare with Van den Dries, Macintyre and Marker in \cite{Dries1994TheET}).
	\item[(3)] A function is log-analytic of order $0$ if and only if it is piecewise the restriction of globally subanalytic functions.
\end{itemize}
\end{rem}

\begin{exc}\label{1.3}
The function 
$$f:\mathbb{R}^2 \to \mathbb{R}, (x,y) \mapsto \arctan(\log(\max\{\log(x^4+\log(y^2+2)),1\})),$$ 
is log-analytic of order (at most) $3$.
\end{exc}
\begin{rem}\label{1.4}
Let $r \in \mathbb{N}_0$. The set of log-analytic functions on $X$ of order at most $r$ as well as the set of log-analytic functions on $X$ is an $\mathbb{R}$-algebra with respect to pointwise addition and multiplication.
\end{rem}
Now we consider a definable function $f$ as a composition of log-analytic functions and exponentials. To formalize the iterations of the exponentials we introduce the \textit{exponential number} of $f$ with respect to a set $E$ of positive definable functions: this is the minimal number of iterations of exponentials from $E$ which are  necessary to write $f$ as a composition of log-analyic functions and functions from $E$ (see also \cite{Dries2002OminimalPT} for such considerations: There the exponential level of a definable function is introduced to describe iterations of exponentials, but without mentioning any set $E$ of positive definable functions.) This allows us to define and study several different classes of definable functions like \textit{$C$-nice functions} (where $E$ consists of $C$-heirs, see Section 2.3) or \textit{restricted log-exp-analytic functions} (where $E$ consists of exponentials of locally bounded functions, see \cite{10.1215/00192082-10234104}). The former helps to understand cellwise preparation of log-analytic functions (see Theorem A) and the latter is a large class of definable functions which contains the log-analytic ones properly and shares some analytic properties with globally subanalytic functions (see \cite{OprisComplex} and \cite{10.1215/00192082-10234104}). %(See also \cite{10.1215/00192082-10234104} for the notion of $C$-consistent and $C$-regular functions which are generalisations of $C$-heirs respectively $C$-nice functions.)

\begin{defn}\label{1.5}
Let $f:X \to \mathbb{R}$ be a function. Let $E$ be a set of positive definable functions on $X$.
\begin{itemize}
	\item [(a)] By induction on $e \in \mathbb{N}_0$ we define that $f$ has \textbf{exponential number at most $e$ with respect to $E$}.
	
	{\bf Base Case}: The function $f$ has exponential number at most $0$ with respect to $E$ if $f$ is log-analytic.
	
	{\bf Inductive Step}: The function $f$ has exponential number at most $e$ with respect to $E$ if the following holds: There are $k,l \in \mathbb{N}_0$, functions $g_1,...,g_k:X \to \mathbb{R}$ and $h_1,...,h_l:X \to \mathbb{R}$ with exponential number at most $e-1$ with respect to $E$ and a log-analytic function $F:\mathbb{R}^{k+l} \to \mathbb{R}$ such that
	$$f=F(g_1,...,g_k,\exp(h_1),...,\exp(h_l))$$
	and $\exp(h_1),...,\exp(h_l) \in E$.
	
	\item [(b)] Let $e \in \mathbb{N}_0$. We say that $f$ has \textbf{exponential number $e$ with respect to $E$} if $f$ has exponential number at most $e$ with respect to $E$ but not at most $e-1$ with respect to $E$.
	
	\item [(c)] We say that $f$ \textbf{can be constructed from $E$} if there is $e \in \mathbb{N}_0$ such that $f$ has exponential number $e$ with respect to $E$. 
\end{itemize}
\end{defn}

\begin{rem}\label{1.6}
Let $E$ be a set of positive definable functions on $X$. Let $f$ be a function on $X$ which can be constructed from $E$. The following holds.
\begin{itemize}
	\item[(1)] $f$ has a unique exponential number $e \in \mathbb{N}_0$ with respect to $E$ and is definable.
	\item[(2)] $f$ has exponential number $0$ with respect to $E$ if and only if $f$ is log-analytic. 
\end{itemize}
\end{rem}

\begin{rem}\label{1.7}
Let $E$ be the set of all positive definable functions on $X$. Then every definable function $f:X \to \mathbb{R}$ can be constructed from $E$.
\end{rem}

\begin{rem}\label{1.8}
Let $e \in \mathbb{N}_0$. Let $E$ be a set of positive definable functions on $X$.
\begin{itemize}
	\item[(1)] Let $f:X \to \mathbb{R}$ be a function with exponential number at most $e$ with respect to $E$. Then $\exp(f)$ has exponential number at most $e+1$ with respect to $E \cup \{\exp(f)\}$.
	\item[(2)] Let $s \in \mathbb{N}_0$. Let $f_1,...,f_s:X \to \mathbb{R}$ be functions with exponential number at most $e$ with respect to $E$ and let $F:\mathbb{R}^s \to \mathbb{R}$ be log-analytic. Then $F(f_1,...,f_s)$ has exponential number at most $e$ with respect to $E$. 
\end{itemize}
\end{rem}

{\bf Proof:}
(1): One sees with Definition \ref{1.5} applied to $g:=F(\exp(f))$ where $F=\id_{\mathbb{R}}$ that $\exp(f)$ has exponential number at most $e+1$ with respect to $E \cup \{\exp(f)\}$.

(2): We may assume $e>0$. Let $k,l \in \mathbb{N}_0$, $g_1,...,g_k,h_1,...,h_l:X \to \mathbb{R}$ be functions with exponential number at most $e-1$ with respect to $E$ with $\exp(h_1),...,\exp(h_l) \in E$, and $G_j:\mathbb{R}^{k+l} \to \mathbb{R}$ be log-analytic such that $f_j=G_j(\beta)$ for $j \in \{1,...,s\}$ where $\beta:=(g_1,...,g_k,\exp(h_1),...,\exp(h_l)).$ Let $v$ range over $\mathbb{R}^{k+l}$. Then
$$H:\mathbb{R}^{k+l} \to \mathbb{R}, v \mapsto F(G_1(v),...,G_s(v)),$$
is log-analytic such that $H(\beta)=F(f_1,...,f_s)$. \hfill$\blacksquare$

\section{A Preparation Theorem for Log-Analytic	Functions}

For the whole section let $n \in \mathbb{N}_0$, $t:=(t_1,...,t_n)$ range over $\mathbb{R}^n$, $x$ over $\mathbb{R}$ and let $\pi:\mathbb{R}^{n+1} \to \mathbb{R}^n, (t,x) \mapsto t$.

\subsection{Logarithmic Scales}
In this subsection we define and investigate \textit{logarithmic scales}. A logarithmic scale is a tuple of functions where each component is defined by taking logarithm of the previous one and then ''translating'' by a definable function depending only on $t$ starting at $C \to C, (t,x) \mapsto x-\Theta_0(t),$ where $\Theta_0$ is definable. This means that the components of a logarithmic scale $\mathcal{Y}(t,-)$ behaves as iterated logarithms independently of $t$. \\
Logarithmic scales are the main technical tools in the cellwise preparation of log-analytic functions and help enormously to describe how log-analytic functions depend on the last variable $x$.

Let $C \subset \mathbb{R}^{n+1}$ be definable. 

%Before introducing the important notion of logarithmic scales, explain it! What is the point here? Definition 2.1 looks very much like the objects introduced in [LR97]. What is the difference? Is is \nicer"? If yes, for which reason?

\begin{defn}\label{2.1}
Let $r \in \mathbb{N}_0$. A tuple $\mathcal{Y}:=(y_0,...,y_r)$ of functions on $C$ is called an \textbf{$r$-logarithmic scale} on $C$ with \textbf{center} $\Theta=(\Theta_0,...,\Theta_r)$ if the following holds:
\begin{itemize}
	\item[(a)] $y_j>0$ or $y_j<0$ for every $j \in \{0,...,r\}$.
	\item[(b)] $\Theta_j$ is a definable function on $\pi(C)$ for every $j \in \{0,...,r\}$.
	\item[(c)] We have $y_0(t,x)=x-\Theta_0(t)$ and inductively $y_j(t,x)=\log(\vert{y_{j-1}(t,x)}\vert) - \Theta_j(t)$ for every $j \in \{1,...,r\}$ and all $(t,x) \in C$.
	\item[(d)] Either there is $\epsilon_0 \in \textnormal{}]0,1[$ such that $0<\vert{y_0(t,x)}\vert < \epsilon_0\vert{x}\vert$ for all $(t,x) \in C$ or $\Theta_0=0$, and for every $j \in \{1,...,r\}$ either there is $\epsilon_j \in \textnormal{}]0,1[$ such that $0<\vert{y_j(t,x)}\vert<\epsilon_j\vert{\log(\vert{y_{j-1}(t,x)}\vert)}\vert$ for all $(t,x) \in C$ or $\Theta_j=0$.
\end{itemize}
\end{defn}
We also write $y_0$ instead of $(y_0)$ for a $0$-logarithmic scale.

Note that Definition 4.1 is a little bit more precise than the objects introduced in \cite{Lion1997}: Lion and Rolin supposed that either $\Theta_j=0$ or there is a positive constant $K$ such that $\vert{y_j(t,x)}\vert \leq K\Theta_j(t)$ for every $(t,x) \in C$ and $j \in \{0,...,n\}$. But for the rigorous poof of the results in \cite{KaiserOpris}, \cite{OprisComplex} and \cite{10.1215/00192082-10234104} the full inequality in (d) is needed. With the objects from \cite{Lion1997} one can for example not prove that $\Theta=0$ in general on a simple cell $C$ (see Proposition 2.19 in \cite{KaiserOpris}, see Definition 2.15 in \cite{KaiserOpris} for the notion of a simple cell $C$).

\begin{rem}\label{2.2}
Let $r \in \mathbb{N}_0$. Let $\mathcal{Y}$ be an \textbf{$r$-logarithmic scale} on $C$ with center $\Theta$. Then $\Theta$ is uniquely determined by $\mathcal{Y}$ and $\mathcal{Y}$ is log-analytic if and only if $\Theta$ is log-analytic.
\end{rem} 

\begin{defn}\label{2.3}
Let $q=(q_0,...,q_r) \in \mathbb{Q}^{r+1}$. We set 
$$\vert{\mathcal{Y}}\vert^{\otimes q}:=\prod_{j=0}^r\vert{y_j}\vert^{q_j}$$
and for $(t,x) \in C$
$$\vert{\mathcal{Y}(t,x)}\vert^{\otimes q}:=\prod_{j=0}^r\vert{y_j(t,x)}\vert^{q_j}.$$
\end{defn}

\begin{defn}\label{2.4}
Let $\mathcal{D}$ be a decomposition of $C$ into finitely many definable cells. Let $f:C \to \mathbb{R}$ be a function. We set
$$C^f:=\{(t,f(t,x)) \in \mathbb{R}^n \times \mathbb{R} \mid (t,x) \in C\}$$
and 
$$\mathcal{D}^f:=\{D^f \mid D \in \mathcal{D}\}.$$	
\end{defn}

\begin{rem}\label{2.5}
The following holds. 
\begin{itemize}
	\item [(1)] Let $f:C \to \mathbb{R}$ be definable. Then $D:=C^f$ is definable. 
	\item [(2)] Let $l \in \{1,...,r\}$. We define $\mu_l:C \to \mathbb{R}, (t,x) \mapsto x-\Theta_l(t),$ and inductively for $j \in \{l+1,...,r\}$ we define $\mu_j:C \to \mathbb{R}, (t,x) \mapsto \log(\vert{\mu_{j-1}(t,x)}\vert)-\Theta_j(t)$. Then $\mathcal{Y}_{r-l,D}:=(\mu_l,...,\mu_r)$ is a well-defined $(r-l)$-logarithmic scale with center $(\Theta_l,...,\Theta_r)$ on $D:=C^{\log(\vert{y_{l-1}}\vert)}$.
\end{itemize} 
\end{rem}

{\bf Proof:}
Property (1) is clear. For property (2) note that
$$\mathcal{Y}(t,x)=(y_0(t,x),...,y_{l-1}(t,x),\mathcal{Y}_{r-l,D}(t,\log(\vert{y_{l-1}(t,x)}\vert)))$$
for every $(t,x) \in C$ where $y_0(t,x)=x -\Theta_0(x)$ and inductively $y_j(t,x)=\log(\vert{y_{j-1}(t,x)}\vert) -\Theta_j(x)$ for $j \in \{1,...,r\}$. So it is straightforward to see with Definition \ref{2.1} that $\mathcal{Y}_{r-l,D}$ is an $(r-l)$-logarithmic scale on $D$. \hfill$\blacksquare$ 

\begin{defn}\label{2.6}
Let $M \in \mathbb{R}_{>0}$. We set
$$C_{>M}(\mathcal{Y}):=\{(t,x) \in C \mid \vert{y_l(t,x)}\vert > M \textnormal{ for all } l \in \{1,...,r\} \},$$
and for every $l \in \{1,...,r\}$
$$C_{l,M}(\mathcal{Y}):=\{(t,x) \in C \mid \vert{y_l(t,x)}\vert \leq M\}.$$
\end{defn}

\begin{rem}\label{2.7}
Let $M \in \mathbb{R}_{>0}$. Then $C_{>M}$ and $C_{1,M},...,C_{r,M}$ are definable.
\end{rem}

\begin{rem}\label{2.8}
{\it Assume $r \in \mathbb{N}$. The following properties hold.
	\begin{itemize}
		\item [(1)] Let $l \in \{1,...,r\}$. We have 
		$$\vert{y_l}\vert \leq \vert{\log(\vert{y_{l-1}}\vert)}\vert$$
		on $C$.
		
		\item [(2)] Let $c,\lambda_1,...,\lambda_r \in \mathbb{R}$. Then there is $M \in \mathbb{R}_{>1}$ such that
		$$\vert{c+\sum_{k=1}^{r}\lambda_k\log(\vert{y_k}\vert)}\vert \leq \frac{\vert{y_1}\vert}{2}$$
		on $C_{>M}$.
	\end{itemize}}
\end{rem}

{\bf Proof:}
Recall that $y_0(t,x)=x -\Theta_0(t)$ and inductively for $j \in \{1,...,r\}$ $y_j(t,x) = \log(\vert{y_{j-1}(t,x)}\vert)-\Theta_j(t)$ for every $(t,x) \in C$.

(1): We find $\epsilon_l \in \textnormal{}]0,1[$ such that $\vert{y_l(t,x)}\vert \leq \epsilon_l\vert{\log(\vert{y_{l-1}(t,x)}\vert)}\vert$ for every $(t,x) \in C$ or $\Theta_l = 0$ by Definition \ref{2.1}. Hence $\vert{y_l}\vert \leq \vert{\log(\vert{y_{l-1}}\vert)}\vert$. 

(2): Take $M \geq \exp_r(1)$ to obtain with $(1)$
$$\vert{c+\sum_{k=1}^r\lambda_k\log(\vert{y_k}\vert)}\vert \leq \vert{c}\vert+\sum_{k=1}^r\vert{\lambda_k}\vert \cdot \log_k(\vert{y_1}\vert)$$
on $C_{>M}$. 
%Zwischenschritt: \vert{c}\vert+\sum_{j=1}^{r}\vert{\lambda_j}\vert \cdot \log(\vert{y_i}\vert) 
By increasing $M$ if necessary we may assume $\vert{c}\vert \leq \frac{\vert{y_1}\vert}{4}$ and $$\vert{\lambda_l}\vert\log_l(\vert{y_1}\vert) \leq \frac{\vert{y_1}\vert}{4r}$$
for every $l \in \{1,...,r\}$ on $C_{>M}$. This gives the result.
\hfill $\blacksquare$

The rest of Subsection 2.2 is important for technical proofs of the main results of this paper.

\begin{defn}\label{2.9}
Let $m \in \mathbb{N}$. Let $D \subset X \subset \mathbb{R}^m$. Let $f,g:X \to \mathbb{R}$ be functions. We call $f$ \textbf{similar to $g$ on $D$}, written $f \sim_D g$, if there is $\delta \in \mathbb{R}_{>0}$ such that $1/\delta \cdot g < f < \delta \cdot g$ on $D$.
\end{defn}

\begin{rem}\label{2.10}
Let $D \subset X \subset \mathbb{R}^m$. Let $f,g:X \to \mathbb{R}$ be functions. If $f \sim_D g$ then $f$ and $g$ don't have a zero on $D$.
\end{rem}

\begin{rem}\label{2.11}
Let $D \subset X \subset \mathbb{R}^m$ and let $f,g:X \to \mathbb{R}$ be functions.
\begin{itemize}
	\item [(1)]
	It holds $f \sim_D g$ if and only if $\frac{f}{g} \sim_D 1$.
	\item [(2)]
	If $f \sim_D g$ then $\textnormal{sign}(f)=\textnormal{sign}(g)$ on $D$.
	\item [(3)]
	The relation $f \sim_D g$ is an equivalence relation on the set of all functions on $X$ without a zero on $D$.
	\item [(4)]
	The set of functions on $X$ which are similar to $1$ on $D$ form a divisible group with respect to pointwise multiplication.
\end{itemize}
\end{rem}

\begin{rem}\label{2.12}
If $\Theta_0 \neq 0$ then $x \sim_C \Theta_0$. Let $j \in \{1,...,r\}$. If $\Theta_j \neq 0$ then $\log(\vert{y_{j-1}}\vert) \sim_C \Theta_j$. 
\end{rem}

{\bf Proof:}
Assume $\Theta_0 \neq 0$. Then there is $\epsilon \in \textnormal{}]0,1[$ such that $\vert{x-\Theta_0(t)}\vert<\epsilon\vert{x}\vert$ for every $(t,x) \in C$. This gives that $x \neq 0$ and
$$1-\epsilon<\frac{\Theta_0(t)}{x}<1+\epsilon$$
for every $(t,x) \in C$. Set $\delta:=\max\{1/(1-\epsilon),1+\epsilon\}$. Then 
$1/\delta<\frac{\Theta_0(t)}{x}<\delta$ for every $(t,x) \in C$.

Assume $\Theta_j \neq 0$. Then there is $\epsilon_j \in \textnormal{}]0,1[$ such that $$\vert{\log(\vert{y_{j-1}(t,x)}\vert)-\Theta_j(t)}\vert<\epsilon_j\vert{\log(\vert{y_{j-1}(t,x)}\vert)}\vert$$
for every $(t,x) \in C$. This gives $\log(\vert{y_{j-1}(t,x)}\vert) \neq 0$ for every $(t,x) \in C$. Now proceed in the same way as above.
\hfill$\blacksquare$

\begin{prop}\label{2.13}
Let $\Psi:\pi(C) \to \mathbb{R}$ be definable. The following properties hold.
\begin{itemize}
		\item [(1)]
		Let $l \in \{1,...,r\}$. If $\vert{y_{l-1}}\vert \sim_{C_{>M}} \Psi$ for some $M >1$ then there is $N \geq M$ such that $\vert{y_l}\vert \sim_{C_{>N}} \vert{\log(\Psi)-\Theta_l}\vert$.
		
		\item [(2)]
		Let $E$ be a set of positive definable functions on $\pi(C)$ such that $\Psi$ and $\Theta_1,...,\Theta_r$ can be constructed from $E$. Let $q_1,...,q_r \in \mathbb{Q}$.
		\begin{itemize}
			\item [(i)] Let $\vert{y_1}\vert \sim_{C_{>M}} \Psi$ for some $M > 1$. There is $N \geq M$ and a function $\mu:\pi(C) \to \mathbb{R}_{>0}$ which can be constructed from $E$ such that 
			$$\prod_{j=1}^r\vert{y_j}\vert^{q_j} \sim_{C_{>N}} \mu.$$
			\item [(ii)] Suppose 
			$$y_0 \sim_C \Psi \prod_{j=1}^r \vert{y_j}\vert^{q_j}.$$
			Then there is $M \in \mathbb{R}_{>1}$ such that $\vert{y_1}\vert \sim_{C_{>M}} \vert{\log(\vert{\Psi}\vert)-\Theta_1}\vert$. Additionally there is $N \geq M$ and a function $\xi:\pi(C) \to \mathbb{R}$ which can be constructed from $E$ such that $y_0 \sim_{C_{>N}} \xi$.
		\end{itemize}
	
\end{itemize}
\end{prop}

{\bf Proof:}
(1): Let $\delta>0$ be with 
$$\frac{1}{\delta} \vert{y_{l-1}}\vert< \Psi <\delta \vert{y_{l-1}}\vert$$
on $C_{>M}$. Note that $\delta>1$. By taking logarithm and subtracting $\Theta_l$ we get
$$-\log(\delta) + y_l < \log(\Psi) -\Theta_l < \log(\delta) + y_l$$
on $C_{>M}$. Set $N:=\max\{M,2\log(\delta)\}$. If $y_l>0$ on $C$ we obtain on $C_{>N}$ that $-\log(\delta)+y_l > 0$ and therefore
$$\vert{y_l}\vert - \log(\delta) < \vert{\log(\Psi) -\Theta_l}\vert < \vert{y_l}\vert + \log(\delta).$$
If $y_l<0$ on $C$ we obtain on $C_{>N}$ that $\log(\delta)+y_l < 0$ and therefore
$$\vert{y_l}\vert - \log(\delta) \leq \vert{\log(\delta) + y_l}\vert < \vert{\log(\Psi) -\Theta_l}\vert < \vert{-\log(\delta)+y_l}\vert \leq \vert{y_l}\vert + \log(\delta).$$
In both cases we obtain
$$\frac{\vert{y_l}\vert}{2} < \vert{\log(\Psi) -\Theta_l}\vert < 2\vert{y_l}\vert$$
on $C_{>N}$.

(2),(i): Let $\Psi_1:=\Psi$. With (1) we find inductively for $l \in \{2,...,r\}$ a real number $N_l \geq N_{l-1}$ such that $$\vert{y_l}\vert \sim_{C_{>N_l}} \vert{\log(\Psi_{l-1}) - \Theta_l}\vert:=\Psi_l$$
where $N_1:=M$. We see with an easy induction on $l \in \{1,...,r\}$ and Remark \ref{1.8} that $\Psi_l$ can be constructed from $E$ for every $l \in\{1,...,r\}$. For $N:=N_r$ we obtain with (4) in Remark \ref{2.11}
$$\vert{y_1}\vert^{q_1} \cdot ... \cdot \vert{y_r}\vert^{q_r} \sim_{C_{>N}} \Psi_1^{q_1} \cdot ... \cdot \Psi_r^{q_r}:=\mu.$$
Note that $\mu:\pi(C) \to \mathbb{R}_{>0}$ can be constructed from $E$ by Remark \ref{1.8}. 

(2),(ii): Let $\delta > 0$ be such that 
$$
\frac{1}{\delta} \vert{y_1}\vert^{q_1} \cdot ... \cdot \vert{y_r}\vert^{q_r} \Psi < y_0 < \delta \vert{y_1}\vert^{q_1} \cdot ... \cdot \vert{y_r}\vert^{q_r} \Psi
$$
on $C$. Let $\kappa:=\max\{1/\delta,\delta\}$. Taking absolute values, logarithm and subtracting $\Theta_1$ we obtain
$$-\log(\kappa) + L + \log(\vert{\Psi}\vert) -\Theta_1 <  y_1  < \log(\kappa) + L + \log(\vert{\Psi}\vert) -\Theta_1$$
on $C$ where $L:=\sum_{k=1}^r q_k\log(\vert{y_k}\vert)$. By Remark \ref{2.8} we find $M>1$ such that 
$$\log(\kappa)+\vert{L}\vert \leq \frac{\vert{y_1}\vert}{2}$$
on $C_{>M}$. We obtain
$$-\frac{\vert{y_1}\vert}{2} + \log(\vert{\Psi}\vert) -\Theta_1 <  y_1  < \frac{\vert{y_1}\vert}{2} + \log(\vert{\Psi}\vert) -\Theta_1$$
and therefore
$$\frac{\vert{y_1}\vert}{2} < \vert{\log(\vert{\Psi}\vert)-\Theta_1}\vert < 2\vert{y_1}\vert$$
on $C_{>M}$. Let $\Gamma:=\vert{\log(\vert{\Psi}\vert)-\Theta_1}\vert$. Then $\Gamma:\pi(C) \to \mathbb{R}_{>0}$ can be constructed from $E$, and it holds $\vert{y_1}\vert \sim_C \Gamma$. By $(2),(i)$ we find a function $\mu:\pi(C) \to \mathbb{R}_{>0}$ which can be constructed from $E$ and $N \geq M$ such that
$$\vert{y_1}\vert^{q_1} \cdot ... \cdot \vert{y_r}\vert^{q_r} \sim_{C_{>N}} \mu.$$
Therefore with (1) in Remark \ref{2.11} $\Psi \cdot \vert{y_1}\vert^{q_1} \cdot ... \cdot \vert{y_r}\vert^{q_r} \sim_{C_{>N}} \Psi \cdot \mu$ and by (3) in Remark \ref{2.11} $y_0 \sim_{C_{>N}} \Psi \cdot \mu.$ So take $\xi:= \Psi \cdot \mu$. Then $\xi$ can be constructed from $E$. \hfill$\blacksquare$

\subsection{A Preparation Theorem for Log-Analytic Functions}

In this section we present the results from \cite{Lion1997}, i.e. preparation theorems for globally subanalytic functions and log-analytic functions. 

\begin{defn}[Lion-Rolin, \cite{Lion1997}, Section 0.4]\label{2.14}
Let $C \subset \mathbb{R}^n \times \mathbb{R}$ be globally subanalytic. Let $f:C \to \mathbb{R}$ be a function. Then $f$ is called \textbf{globally subanalytically prepared in $x$ with center $\theta$} if for every $(t,x) \in C$ 
$$f(t,x)=a(t)\cdot \vert{x-\theta(t)}\vert^q \cdot u(t,x)$$
where $q \in \mathbb{Q}$, $\theta:\pi(C) \to \mathbb{R}$ is a globally subanalytic function such that either $x > \theta(t)$ for every $(t,x) \in C$ or $x < \theta(t)$ for every $(t,x) \in C$, $a$ is a globally subanalytic function on $\pi(C)$ which is identically zero or has no zeros, and the following holds for $u$. There is $s \in \mathbb{N}$ such that $u=v \circ \phi$ where $v$ is a power series which converges on an open neighborhood of $[-1,1]^s$ with $v([-1,1]^s) \subset \mathbb{R}_{>0}$ and $\phi:=(\phi_1,...,\phi_s): C \to [-1,1]^s$ with $\phi_j(t,x)=b_j(t)\vert{x-\theta(t)}\vert^{p_j}$ for $(t,x) \in C$ where $p_j \in \mathbb{Q}$ and $b_j:\pi(C) \to \mathbb{R}$ is globally subanalytic for $j \in \{1,...,s\}$. Additionally either there is $\epsilon \in \textnormal{}]0,1[$ such that 
$$\vert{x-\theta(t)}\vert < \epsilon \vert{x}\vert$$
for every $(t,x) \in C$ or $\theta=0$.
\end{defn}

We see that a globally subanalytic prepared function has roughly speaking the form of a Puiseux series in one variable. 

\begin{fac}[Lion-Rolin, \cite{Lion1997}, Section 1]\label{2.15} 
Let $m \in \mathbb{N}$. Let $X \subset \mathbb{R}^n \times \mathbb{R}$ and $f_1,...,f_m:X \to \mathbb{R}$ be globally subanalytic. Then there is a globally subanalytic cell decomposition $\mathcal{C}$ of $X_{\neq 0}$ such that for every $C \in \mathcal{C}$ there is a globally subanalytic $\theta:\pi(C) \to \mathbb{R}$ such that $f_1,...,f_m$ are globally subanalytically prepared in $x$ with center $\theta$ on $C$.
\end{fac}
\begin{rem}\label{2.16}
There are similar versions of this preparation theorem for reducts of $\mathbb{R}_{\textnormal{an}}$ like $\mathbb{R}_{\mathcal{W}}$ where $\mathcal{W}$ is a convergent Weierstrass system (compare with \cite{Weiherstrass} for the details).
\end{rem}
Let $C \subset \mathbb{R}^n \times \mathbb{R}$ be a definable cell.

\begin{defn}[Lion-Rolin, \cite{Lion1997}, Section 0.4]\label{2.17}
Let $r \in \mathbb{N}_0$. Let $g:C \to \mathbb{R}$ be a function. We say that $g$ is \textbf{$r$-log-analytically prepared in $x$ with center $\Theta$} if
$$g(t,x)=a(t) \vert{\mathcal{Y}(t,x)}\vert^{\otimes q}u(t,x)$$
for all $(t,x) \in C$ where $a$ is a definable function on $\pi(C)$ which vanishes identically or has no zero, $\mathcal{Y}=(y_0,...,y_r)$ is an $r$-logarithmic scale with center $\Theta$ on $C$, $q \in \mathbb{Q}^{r+1}$ and the following holds for $u$. There is $s \in \mathbb{N}$ such that $u=v \circ \phi$ where $v$ is a power series which converges on an open neighborhood of $[-1,1]^s$ with $v([-1,1]^s) \subset \mathbb{R}_{>0}$ and $\phi:=(\phi_1,...,\phi_s):C \to [-1,1]^s$ is a function of the form 
$$\phi_j(t,x):=b_j(t)\vert{\mathcal{Y}(t,x)}\vert^{\otimes p_j}$$
for $j \in \{1,...,s\}$ and $(t,x) \in C$ where $b_j:\pi(C) \to \mathbb{R}$ is definable and $p_j:=(p_{j0},...,p_{jr}) \in \mathbb{Q}^{r+1}$. We call $a$ \textbf{coefficient} and $b:=(b_1,...,b_s)$ a tuple of \textbf{base functions} for $f$. An \textbf{LA-preparing tuple} for $f$ is then
$$\mathcal{J}:=(r,\mathcal{Y},a,q,s,v,b,P)$$
where
$$P:=\left(\begin{array}{cccc}
p_{10}&\cdot&\cdot&p_{1r}\\
\cdot&& &\cdot\\
\cdot&& &\cdot\\
p_{s0}&\cdot&\cdot&p_{sr}\\
\end{array}\right)\in M\big(s\times (r+1),\mathbb{Q}).$$
\end{defn}

\begin{rem}\label{2.18}
Let $f:C \to \mathbb{R}$ be a function. Let $f\sim_C g$ where $g$ is $r$-log-analytically prepared in $x$ with preparing tuple $(r,\mathcal{Y},a,q,s,v,b,P)$. Then
$$f \sim_C a \vert{\mathcal{Y}}\vert^{\otimes q}.$$
\end{rem}

\begin{defn}\label{2.19}
Let $r \in \mathbb{N}_0$ and $m \in \mathbb{N}$. Let $g_1,...,g_m:C \to \mathbb{R}$ be functions. We say that $g_1,...,g_m$ are $r$-log-analytically prepared in $x$ \textbf{in a simultaneous way} if there is $\Theta:\pi(C) \to \mathbb{R}^{r+1}$ such that $g_1,...,g_m$ are $r$-log-analytically prepared in $x$ with center $\Theta$.
\end{defn}

\begin{rem}\label{2.20}
Let $m \in \mathbb{N}$ and $r \in \mathbb{N}_0$. Let $f_1,...,f_m:C \to \mathbb{R}$ be $r$-log-analytically prepared in $x$ in a simultaneous way. Then there are preparing tuples for $f_1,...,f_m$ which coincide in $r,\mathcal{Y},s,b,P$.
\end{rem}

{\bf Proof:}
This follows immediately with Definition \ref{2.17} and by redefining the corresponding power series $v_1,...,v_m$. \hfill$\blacksquare$

\begin{fac}[Lion-Rolin, \cite{Lion1997}, Section 2.1]\label{2.21}
Let $X \subset \mathbb{R}^n \times \mathbb{R}$ be definable. Let $r \in \mathbb{N}_0$ and $m \in \mathbb{N}$. Let $f_1,...,f_m:X \to \mathbb{R}$ be log-analytic functions of order at most $r$. Then there is a definable cell decomposition $\mathcal{C}$ of $X_{\neq 0}$ such that $f_1|_C,...,f_m|_C$ are $r$-log-analytically prepared in $x$ in a simultaneous way for every $C \in \mathcal{C}$.
\end{fac}

The rest of this section is about the following question: Are the coefficient, center and base functions of the cellwise preparation of a log-analytic function again log-analytic? We will see that this question generally has a negative answer.

\begin{defn}\label{2.22}
Let $\mathcal{Y}$ be an $r$-logarithmic scale on $C$ with center $\Theta=(\Theta_0,...,\Theta_r)$. We call $\mathcal{Y}$ \textbf{pure} if its center $\Theta$ is log-analytic.
\end{defn}

\begin{defn}\label{2.23}
Let $r \in \mathbb{N}_0$.

\begin{itemize}
\item [(a)] Let $C$ be a cell and $f:C \to \mathbb{R}$ be a function. We call $f:C \to \mathbb{R}$ \textbf{purely $r$-log-analytically prepared in $x$ with center $\Theta$} if $f$ is $r$-log-analytically prepared in $x$ with log-analytic center $\Theta$, log-analytic coefficient and log-analytic base functions. An LA-preparing tuple for $g$ with log-analytic components is then called a \textbf{pure LA-preparing tuple} for $g$. 

\item [(b)] Let $f_1,...,f_m:C \to \mathbb{R}$ be functions. We call $f_1,...,f_m$ purely $r$-log-analytically prepared in $x$ \textbf{in a simultaneous way} if $f_1,...,f_m$ are purely $r$-log-analytically prepared in $x$ with the same center.
\end{itemize}  
\end{defn}

\begin{rem}\label{2.24}
Let $r \in \mathbb{N}_0$. If $g$ is $r$-log-analytically prepared in $x$ then $g$ is definable but not necessarily log-analytic. If $g$ is purely $r$-log-analytically prepared in $x$ then $g$ is log-analytic.
\end{rem}

Purely $r$-log-analytically prepared functions in $x$ induce the following definition of \textit{purely log-analytic functions in $x$}.

\begin{defn}\label{2.25}
Let $X \subset \mathbb{R}^n \times \mathbb{R}$ be definable and let $f:X \to \mathbb{R}$ be a function in $x$.
\begin{itemize}
	\item[(a)] Let $r \in \mathbb{N}_0$. By induction on $r$ we define that $f$ is \textbf{purely log-analytic in $x$ of order at most $r$}.
	
	\textbf{Base case}: The function $f$ is purely log-analytic in $x$ of order at most $0$ if the following holds: There is a definable cell decomposition $\mathcal{C}$ of $X$ such that for every $C \in \mathcal{C}$ there is $m \in \mathbb{N}_0$, a globally subanalytic $F:\mathbb{R}^{m+1} \to \mathbb{R}$ and a log-analytic $g:\pi(C) \to \mathbb{R}^m$ such that $f(t,x)=F(g(t),x)$ for $(t,x) \in C$.
	
	\textbf{Inductive step}: The function $f$ is purely log-analytic in $x$ of order at most $r$ if the following holds: There is a definable cell decomposition $\mathcal{C}$ of $X$ such that for $C \in \mathcal{C}$ there are $k,l \in \mathbb{N}_0$, purely log-analytic functions $g_1,...,g_k:C \to \mathbb{R}, h_1,...,h_l:C \to \mathbb{R}_{>0}$ in $x$ of order at most $r-1$, and a globally subanalytic function $F:\mathbb{R}^{k+l} \to \mathbb{R}$ such that
	$$f|_C=F(g_1,...,g_k,\log(h_1),...,\log(h_l)).$$
	\item[(b)] Let $r \in \mathbb{N}_0$. The function $f$ is \textbf{purely log-analytic in $x$ of order $r$} if $f$ is purely log-analytic in $x$ of order at most $r$ but not purely log-analytic in $x$ of order at most $r-1$.
	\item[(c)] The function $f$ is \textbf{purely log-analytic in $x$} if there is $r \in \mathbb{N}_0$ such that $f$ is purely log-analytic in $x$ of order $r$.
\end{itemize}
\end{defn}

\begin{rem}\label{2.26}
Let $X \subset \mathbb{R}^n \times \mathbb{R}$ be definable. A purely log-analytic function in $x$ on $X$ is log-analytic.
\end{rem}

\begin{rem}\label{2.27}
Let $f:C \to \mathbb{R}$ be a function. Let $r \in \mathbb{N}_0$.
\begin{itemize}
	\item [(1)]
	A log-analytic function $g:C \to \mathbb{R}$ of order at most $r$ is purely log-analytic in $x$ of order at most $r$.
	\item [(2)]
	A purely $r$-log-analytically prepared function $g:C \to \mathbb{R}$ in $x$ is purely log-analytic in $x$ of order at most $r$.
\end{itemize}
\end{rem}

Now we can give the promised example that the above Fact \ref{2.21} can in general not be carried out in the log-analytic category.

\begin{exc}[Kaiser - Opris, \cite{KaiserOpris}]\label{2.28}
Let $\phi: \textnormal{}]0,\infty[\to \mathbb{R}, y \mapsto y/(1+y)$. Consider the log-analytic function 
$$f:\mathbb{R}_{>0} \times \mathbb{R}, (t,x) \mapsto -\frac{1}{\log(\phi(x))} - t.$$
Then $f$ is log-analytic of order $1$, but does not allow piecewise a pure $1$-log-analytic preparation in $x$. 
\end{exc}

{\bf Proof:}
For the reader's convenience we present the proof from \cite{KaiserOpris} here. Assume that the contrary holds. Let $\mathcal{C}$ be the corresponding cell decomposition. Let $\psi:\textnormal{}]0,1[ \to \mathbb{R}, y \mapsto y/(1-y)$. Then $\psi$ is the compositional inverse of $\phi$. Note that $f(t,\psi(e^{-1/t}))=0$ for all $t \in \mathbb{R}_{>0}$. Let $\alpha:\mathbb{R}_{>0} \to \mathbb{R}, t \mapsto \psi(e^{-1/t})$. Then $\alpha$ is not log-analytic and $\alpha(t)=\sum_{n=1}^{\infty}e^{-n/t}$ for all $t \in \mathbb{R}_{>0}$. By passing to a finer definable cell decomposition we find a cell $C$ of the form 
$$C:=\{(t,x) \in \mathbb{R}_{>0} \times \mathbb{R}_{>0} \mid 0<t<\epsilon, \alpha(t)<x<\alpha(t)+\eta(t)\}$$   
with some suitable $\epsilon \in \mathbb{R}_{>0}$ and some definable function $\eta:\textnormal{}]0,\epsilon[\textnormal{} \to \mathbb{R}_{>0}$ such that $f$ is purely $1$-log-analytically prepared in $x$ on $C$. Let $(1,\mathcal{Y},a,q,s,v,b,P)$ be a pure preparing tuple for $f|_C$ and let $\Theta=(\Theta_0,\Theta_1)$ be the center of $\mathcal{Y}$. 

{\bf Claim}\\
$\Theta_0=0$.

{\bf Proof of the claim:} Assume that $\Theta_0$ is not the zero function. By the definition of a $1$-logarithmic scale we find $\varepsilon_0 \in \textnormal{}]0,1[$ such that $\vert{y_0}\vert<\varepsilon_0 \vert{x}\vert$ on $C$. This implies $\vert{\alpha(t)-\Theta_0(t)}\vert \leq \epsilon_0\alpha(t)$ for all $0<t<\epsilon$. But this is not possible since we have $\lim_{t \searrow 0}\alpha(t)/\Theta_0(t)=0$ by the assumption that $\Theta_0$ is log-analytic and not the zero function.

$\hfill \blacksquare_{\textnormal{Claim}}$

From 
$$f(t,x)=a(t)\vert{\mathcal{Y}(t,x)}\vert^{\otimes q} u(t,x)$$
for all $(t,x) \in C$ (where $u:C \to \mathbb{R}_{>0}$ is suitable with $u \sim_C 1$) and 
$$\lim_{x \searrow \alpha(t)} f(t,x)=0$$
for all $t \in \textnormal{}]0,\epsilon[$ we get by o-minimality that there is, after shrinking $\epsilon>0$ if necessary, some $j \in \{0,1\}$ such that $\lim_{x \searrow \alpha(t)} \vert{y_j(t,x)}\vert^{q_j}=0$ for all $t \in \textnormal{}]0,\varepsilon[$. By the claim the case $j=0$ is not possible. In the case $j=1$ we have, again by the claim, that $q_1>0$ and therefore $\Theta_1=\log(\alpha)$. But this is a contradiction to the assumption that the function $\Theta_1$ is log-analytic, since the function $\log(\alpha)$ on the right is not log-analytic. This can be seen by applying the logarithmic series. We obtain $\lim_{t \searrow 0} (\log(\alpha(t)) + 1/t)/e^{-1/t} \in \mathbb{R} \setminus \{0\}$. \hfill$\blacksquare$

Let $f$, $C$, $\alpha$ and $\psi$ be as in Example \ref{2.28}. In the following we mention that $f$ can be $0$-log-analytically prepared on $C$ (after shrinking $\epsilon$ and $\eta$ if necessary) with coefficient, center and base functions which can be constructed from a set $E$ of positive definable functions with the following property: Every $g \in \log(E)$ is a component of the center of a logarithmic scale on $C$.

For $(w_2,w_3,w_4) \in \mathbb{R}^3$ with $w_2,w_4>0$, $w_4/w_2 \in \textnormal{} ]1/2,3/2[$ and $1/(1+w_4) \in \textnormal{} ]1/2,3/2[$ we set
$$u(w_2,w_3,w_4):=\log^*(\frac{w_4}{w_2})+\log^*(\frac{1}{1+w_4}) + w_3$$
where
$$\log^*:\mathbb{R} \to \mathbb{R}, y \mapsto \left\{\begin{array}{ll} \log(y), & y \in [1/2,3/2], \\
0, & \textnormal{ otherwise.} \end{array}\right.$$

Consider
$$F:\mathbb{R}^4 \to \mathbb{R}, (w_1,...,w_4) \mapsto$$
$$\left\{\begin{array}{ll} -\frac{1}{u(w_2,w_3,w_4)}-w_1, & w_2,w_4>0, \textnormal{ } w_4/w_2 \in \textnormal{} ]1/2,3/2[, \textnormal{ } 1/(1+w_4) \in \textnormal{} ]1/2,3/2[  \\ & \textnormal{and } u(w_2,w_3,w_4) \neq 0, \\
 & \\
0, & \textnormal{otherwise.} \end{array}\right.$$

Note that $F$ is globally subanalytic since $\log^*$ is globally subanalytic. By shrinking $\epsilon$ and $\eta$ if necessary we obtain $x/\alpha(t) \in \textnormal{} ]1/2,3/2[$, $1/(1+x) \in \textnormal{}]1/2,3/2[$ and therefore
$$f(t,x)=F(t,\alpha(t),\log(\alpha(t)),x)$$
for $(t,x) \in C$. Note that $\alpha$ can be constructed from $E:=\{\pi(C) \to \mathbb{R}_{>0}, t \mapsto e^{-1/t}\}$ (since $\alpha(t)=\psi(e^{-1/t})$ for $t \in \pi(C)$ and $\psi$ is globally subanalytic). Therefore $\log(\alpha)$ can also be constructed from $E$ (since $\pi(C) \to \mathbb{R}, t \mapsto \log(\psi(t))$, is log-analytic). With Fact 2.15 we find a globally subanalytic cell decomposition $\mathcal{D}$ of $\mathbb{R}^4$ such that for every $D \in \mathcal{D}$ we have that $F|_D$ is globally subanalytically prepared in $x_4$. By shrinking $\epsilon$ and $\eta$ if necessary we may assume that there is $D \in \mathcal{D}$ such that $(t,\alpha(t),\log(\alpha(t)),x) \in D$ for every $(t,x) \in C$. So we may assume that $F$ is globally subanalytically prepared in $w_4$. Consequently $f$ is $0$-log-analytically prepared in $x$ with coefficient, center and base functions which can be constructed from $E$. By further shrinking $\epsilon$ and $\eta$ if necessary we may assume that 
$$\vert{\log(x)+1/t}\vert < 1/2 \vert{\log(x)}\vert$$
for $(t,x) \in C$. Let $y_0:=x$ and $y_1:=\log(y_0)+1/t$. We have that $(y_0,y_1)$ is a $1$-logarithmic scale on $C$ with center $(\Theta_0,\Theta_1)$ where $\Theta_0:=0$ and $\Theta_1:\pi(C) \to \mathbb{R}, t \mapsto -1/t$. So $E$ is a set of one positive definable function on $\pi(C)$ whose logarithm coincides with $\Theta_1$. 

We will see in Section 2.3 that functions which can be constructed from a set of positive definable functions whose logarithm coincide with a component of the center of a logarithmic scale play a crucial role in preparing log-analytic functions.

\subsection{A Type of Definable Functions Closed Under Log-Analytic Preparation}

The main goal for this section is to prove Theorem A. At first we introduce the notion of \textit{$C$-heirs} and \textit{$C$-nice} functions for a definable cell $C \subset \mathbb{R}^n \times \mathbb{R}$, give elementary properties and examples. A $C$-heir is a positive definable function on $\pi(C)$ whose logarithm coincides with the component of a center of a logarithmic scale on $C$. A $C$-nice function is a definable function which is a composition of log-analytic functions and $C$-heirs (i.e. can be constructed from a set $E$ of $C$-heirs, see (c) in Definition \ref{1.5}). As Theorem A suggests, $C$-nice functions are crucial to understand log-analytic functions from the viewpoint of their cellwise preparation.

For Section 2.3 we need the following additional notation. ''Log-analytically prepared'' means always ''log-analytically prepared in $x$''. For $k,l \in \mathbb{N}$ and $s \in \{1,...,l\}$ denote by $M_s(k \times l,\mathbb{Q})$ the set of all $k \times l$-matrices with rational entries where the first $s$ columns are zero.
%Einführung C-Erbe, C-nett informell erklären!

Let $C \subset \mathbb{R}^n \times \mathbb{R}$ be definable. 

\begin{defn}\label{2.29}
We call $g:\pi(C) \to \mathbb{R}$ \textbf{$C$-heir} if there is $r \in \mathbb{N}_0$, an $r$-logarithmic scale $\mathcal{Y}$ with center $(\Theta_0,...,\Theta_r)$ on $C$, and $l \in \{1,...,r\}$ such that $g=\exp(\Theta_l)$.
\end{defn}

\begin{rem}\label{2.30}
Let $D \subset \mathbb{R}^n \times \mathbb{R}$ be definable with $\pi(C)=\pi(D)$. A $C$-heir is not necessarily a $D$-heir.
\end{rem}
{\bf Proof:}
Let $C:=\mathbb{R}^n \times \textnormal{} ]0,1[$ and $D:=\mathbb{R}^n \times \mathbb{R}_{\neq 0}$. Note that there is no $r$-logarithmic scale on $D$ for every $r \in \mathbb{N}_0$. So a $D$-heir does not exist. But it is straightforward to see that $(y_0,y_1)$ with $y_0:=x$ and $y_1:=\log(x)$ is a $1$-logarithmic scale with center $0$ on $C$ and therefore that $g:\pi(C) \to \mathbb{R},x \mapsto 1,$ is a $C$-heir. \hfill$\blacksquare$

\begin{defn}\label{2.31}
We call $g:\pi(C) \to \mathbb{R}$ \textbf{$C$-nice} if there is a set $E$ of $C$-heirs such that $g$ can be constructed from $E$.
\end{defn}

\begin{exc}\label{2.32}
The following holds:
\begin{itemize}
	\item [(1)] A log-analytic function $f:\pi(C) \to \mathbb{R}$ is $C$-nice.
	\item [(2)] Let $C:= \textnormal{}]0,1[^2$ and let $h: \textnormal{}]0,1[ \to \mathbb{R}, t \mapsto e^{-1/t}$. Then $h$ is not $C$-nice.
\end{itemize}
\end{exc}

{\bf Proof:}
(1): Note that $f$ can be constructed from $E=\emptyset$. Therefore $f$ is $C$-nice by Definition \ref{2.31}.\\
(2): Note that $\pi(C)= \textnormal{}]0,1[$. Suppose that $h$ is $C$-nice. Let $E$ be a set of $C$-heirs such that $h$ can be constructed from $E$. An easy calculation shows that the center of every logarithmic scale on $C$ vanishes (compare with Definition 2.1(d) or Proposition 2.19 in \cite{KaiserOpris}). So we see $E=\emptyset$ or $E=1$.

{\bf Claim}\\
The function $h$ is log-analytic.

{\bf Proof of the claim:} Suppose $E = \{1\}$. Let $e \in \mathbb{N}_0$ be such hat $h$ has exponential number at most $e$ with respect to $E$. We do an induction on $e$. For $e=0$ this is clear with part (a) in Definition \ref{1.5}. So suppose $e>0$. Let $u_1,...,u_k,v_1,...,v_l:\pi(C) \to \mathbb{R}$ be functions with $\exp(v_1),...,\exp(v_l) \in E$ which have exponential number at most $e-1$ with respect to $E$ and $F:\mathbb{R}^k \times \mathbb{R}^l \to \mathbb{R}$ be log-analytic such that  
$$h(t)=F(u_1(t),...,u_k(t),\exp(v_1(t)),...,\exp(v_l(t)))$$
for $t \in \pi(C)$. Note that $u_1,...,u_k$ are log-analytic by the inductive hypothesis and that $\exp(v_1)=...=\exp(v_l)=1$. So we obtain that $h$ is log-analytic. If $E=\emptyset$ one sees immediately with part (a) in Definition \ref{1.5} that $h$ is log-analytic as a consequence that that $h$ is $C$-nice. \hfill$\blacksquare_{\textnormal{Claim }}$

But we have $\lim_{t \to 0} h(t)/e^{-1/t} \in \mathbb{R} \setminus \{0\}$, a contradiction to the claim since every log-analytic function is polynomially bounded. \hfill$\blacksquare$

A $C$-nice function which is not log-analytic can be found in Example \ref{2.37}.

\begin{rem}\label{2.33}
Let $r \in \mathbb{N}_0$. Let $(\Theta_0,...,\Theta_r)$ be a $C$-nice center of an $r$-logarithmic scale $\mathcal{Y}$ on $C$. Let $f=\exp(\Theta_j)$ for $j \in \{1,...,r\}$. Then $f$ is a $C$-nice $C$-heir.  
\end{rem}

{\bf Proof:}
Let $E$ be a set of $C$-heirs such that $\Theta_1,...,\Theta_r$ can be constructed from $E$. Then $\exp(\Theta_j)$ can be constructed from $E \cup \{\exp(\Theta_j)\}$ for $j \in \{1,...,r\}$ by Remark \ref{1.8} and is therefore $C$-nice, because $\pi(C) \mapsto \mathbb{R}, t \mapsto \exp(\Theta_j(t)),$ is a $C$-heir. \hfill$\blacksquare$

\begin{rem}\label{2.34}
Let $D \subset C$ be definable. 
\begin{itemize}
	\item [(1)] Let $g:\pi(C) \to \mathbb{R}$ be a $C$-heir. Then $g|_{\pi(D)}$ is a $D$-heir.
	\item [(2)] Let $h:\pi(C) \to \mathbb{R}$ be $C$-nice. Then $h|_{\pi(D)}$ is $D$-nice.
\end{itemize} 
\end{rem}

{\bf Proof:}
(1): This follows from the following fact: Let $r \in \mathbb{N}_0$. Let $\mathcal{Y}$ be an $r$-logarithmic scale on $C$ with center $(\Theta_0,...,\Theta_r)$. Then  $\mathcal{Y}|_D$ is an $r$-logarithmic scale with center $(\Theta_0|_{\pi(D)},...,\Theta_r|_{\pi(D)})$ on $D$.\\
(2): Let $E$ be a set of $C$-heirs such that $h$ can be constructed from $E$. Then it is easily seen that $h|_{\pi(D)}$ can be constructed from $E|_{\pi(D)}$ which is a set of $\pi(D)$-heirs by (1). \hfill$\blacksquare$

\begin{rem}\label{2.35}
The set of $C$-nice functions is closed under composition with log-analytic functions, i.e. let $m \in \mathbb{N}$ and $F:\mathbb{R}^m \to \mathbb{R}$ be log-analytic and $\eta:=(\eta_1,...,\eta_m):\pi(C) \to \mathbb{R}^m$ be $C$-nice. Then $F \circ \eta:\pi(C) \to \mathbb{R}$ is $C$-nice.
\end{rem}

{\bf Proof:}
Let $E$ be a set of $C$-heirs such that $\eta_1,...,\eta_m$ can be constructed from $E$. By Remark \ref{1.8} $F(\eta_1,...,\eta_m)$ can be constructed from $E$.   
\hfill$\blacksquare$

\begin{defn}\label{2.36}
Let $r \in \mathbb{N}_0$. Let $\mathcal{Y}$ be an $r$-logarithmic scale on $C$. Then $\mathcal{Y}$ is \textbf{nice} if its center $\Theta:=(\Theta_0,...,\Theta_r)$ is $C$-nice.
\end{defn}

The next example shows that not every nice logarithmic scale is log-analytic in general. 

\begin{exc}\label{2.37}
Consider the definable cell
$$C:=\{(t,x) \in \mathbb{R}^2 \mid 0<t<1, \textnormal{} \tfrac{1}{1+t}+e^{-2/t+2e^{-1/t}}<x<\tfrac{1}{1+t}+e^{-1/t}\}.$$
Then there is a nice logarithmic scale on $C$ which is not log-analytic. 
\end{exc}

{\bf Proof:}
We have $\pi(C)=\textnormal{}]0,1[$, because
$$e^{-2/t+2e^{-1/t}}<e^{-1/t}$$
for every $t \in \textnormal{}]0,1[$. Consider
$$\Theta_0:\pi(C) \to \mathbb{R}, t \mapsto \tfrac{1}{1+t},$$
$$\Theta_1:\pi(C) \to \mathbb{R}, t \mapsto -1/t,$$
and
$$\hat{\Theta}_1:\pi(C) \to \mathbb{R}, t \mapsto -1/t+e^{-1/t}.$$
For $(t,x) \in C$ consider 
$$y_0(t,x):=x-\Theta_0(t),$$
$$y_1(t,x):=\log(x-\Theta_0(t))-\Theta_1(t),$$
and
$$\hat{y}_1(t,x):=\log(x-\Theta_0(t))-\hat{\Theta}_1(t).$$

{\bf Claim}\\
$\mathcal{Y}:=(y_0,y_1)$ and $\hat{\mathcal{Y}}:=(y_0,\hat{y}_1)$ are $1$-logarithmic scales on $C$. 

{\bf Proof of the claim:} We have $y_0>0$, $y_1<0$, and $\hat{y}_1<0$ on $C$. Let $\epsilon_0:=1/2$ and $\epsilon_1:=1/2$. Let $(t,x) \in C$. 
We have $x < \frac{2}{1+t}$, because $e^{-1/t}<\frac{1}{1+t}$. Therefore $\vert{x-\Theta_0(t)}\vert < \epsilon_0 \vert{x}\vert$. Note that 
$$e^{-2/t+2e^{-1/t}}+\frac{1}{1+t}<x.$$
So an easy calculation shows that
$$\vert{\log(x-\Theta_0(t))-\Theta_1(t)}\vert < \epsilon_1 \cdot \vert{\log(x - \Theta_0(t))}\vert$$
and
$$\vert{\log(x-\Theta_0(t))-\hat{\Theta}_1(t)}\vert < \epsilon_1 \cdot \vert{\log(x - \Theta_0(t))}\vert.$$
\hfill$\blacksquare_{\textnormal{Claim}}$
\iffalse
$$-\log(y-\Theta_0(x))+\Theta_1(x) < 0,5 \cdot \vert{\log(\vert{y - \Theta_0(x)}\vert)}\vert$$
$$\Leftrightarrow -\log(y-\Theta_0(x))+\Theta_1(x) < -0,5 \log (y-\Theta_0(x))$$
$$\Leftrightarrow \Theta_1(x) < 0,5 \cdot \log(y-\Theta_0(x))$$
$$\Leftrightarrow e^{2\Theta_1(x)}+\Theta_0(x) < y$$
$$\Leftrightarrow e^{-2/x} + \tfrac{1}{1+x} < y.$$
\fi \iffalse
$$-\log(y-\Theta_0(x))+\hat{\Theta}_1(x) < 0,5 \cdot \vert{\log(\vert{y - \Theta_0(x)}\vert)}\vert$$
$$\Leftrightarrow e^{2\hat{\Theta}_1(x)}+\Theta_0(x) < y$$
$$\Leftrightarrow e^{-2/x+2e^{-1/x}}+\tfrac{1}{1+x} < y$$
\fi

Note that $\Theta_0$ and $\Theta_1$ are log-analytic. Therefore $\mathcal{Y}$ is nice by (1) in Example \ref{2.32}. Note also that $\hat{\Theta}_1$ is not log-analytic. We have $\hat{\Theta}_1=G(1/t,e^{-1/t})$ where $G:\mathbb{R}^2 \to \mathbb{R}, (w_1,w_2) \mapsto -w_1+w_2,$ is log-analytic. So $\hat{\Theta}_1$ has exponential number $1$ with respect to $E:=\{g\}$ where $g:\pi(C) \to \mathbb{R}, t \mapsto e^{-1/t}$  (since $t \mapsto -1/t$ is globally subanalytic). So $\hat{\Theta}_1$ can be constructed from $E$. By the claim we see that $g$ is a $C$-heir. Therefore $\hat{\mathcal{Y}}$ is an example for a nice $1$-logarithmic scale which is not log-analytic. \hfill$\blacksquare$

\begin{rem}\label{2.38}
Let $r \in \mathbb{N}_0$. Let $\mathcal{Y}:=(y_0,...,y_r)$ be a $C$-nice $r$-logarithmic scale with center $(\Theta_0,...,\Theta_r)$ on $C$. Let $\Psi:\pi(C) \to \mathbb{R}$ be $C$-nice. Assume that
$$y_0 \sim_C \Psi \prod_{j=1}^r \vert{y_j}\vert^{q_j}$$
where $q_1,...,q_r \in \mathbb{Q}$. Then there is $M \in \mathbb{R}_{>1}$ and a $C$-nice $\xi:\pi(C) \to \mathbb{R}$ such that $y_0 \sim_{C_{>M}} \xi$.
\end{rem}

{\bf Proof:}
Let $E$ be a set of $C$-heirs such that $\Psi,\Theta_0,...,\Theta_r$ can be constructed from $E$. The assertion follows from (2),(ii) in Proposition \ref{2.13}. \hfill$\blacksquare$

Now we are able to introduce the notion of nice log-analytic prepared functions.

\begin{defn}\label{2.39}
Let $f:C \to \mathbb{R}$ be a function. We say that $f$ is \textbf{nicely $r$-log-analytically prepared with center $\Theta$} if $f$ is $r$-log-analytically prepared with a nice $r$-logarithmic scale $\mathcal{Y}$ with center $\Theta$, $C$-nice coefficient and $C$-nice base functions. A corresponding LA-preparing tuple for $f$ is then called a \textbf{nice LA-preparing tuple} for $f$.
\end{defn}
\begin{rem}\label{2.40}
Let $r,m \in \mathbb{N}$ and $k \in \{1,...,r\}$. Let $\mathcal{Y}_{k-1}:=(y_0,...,y_{k-1})$ be a nice $(k-1)$-logarithmic scale with center $(\Theta_0,...,\Theta_{k-1})$ on $C$. Let $B:=C^{\log(\vert{y_{k-1}}\vert)}$. Let $\alpha_1,...,\alpha_m:B \to \mathbb{R}$ be nicely $(r-k)$-log-analytically prepared with center $(\Theta_k,...,\Theta_r)$. For every $j \in \{1,...,m\}$ consider
$$\beta_j: C \to \mathbb{R}, (t,x) \mapsto \alpha_j(t,\log(\vert{y_{k-1}(t,x)}\vert)).$$
Let $\Theta:=(\Theta_0,...,\Theta_r)$. Then $\beta_1,...,\beta_m$ are nicely $r$-log-analytically prepared with center $\Theta$. Additionally there is a nice LA-preparing tuple $(r,\mathcal{Y}_r,a_j,q_j,s,v_j,b,P)$ for $\beta_j$ such that $q_j \in \{0\}^k \times \mathbb{Q}^{r+1-k}$ and $P \in M_k(s \times (r+1),\mathbb{Q})$ where $\mathcal{Y}_r$ is the $r$-logarithmic scale on $C$ with center $\Theta:=(\Theta_0,...,\Theta_r)$ for $j \in \{1,...,m\}$.
\end{rem}

{\bf Proof:}
Let 
$$(r-k,\hat{\mathcal{Y}}_{r-k,B},a_j,q_j,s,v_j,b,P)$$
be a nice LA-preparing tuple for $\alpha_j$ where $j \in \{1,...,m\}$. (By Remark \ref{2.20} $s,b,P$ are independent from $j$.) Here $\hat{\mathcal{Y}}_{r-k,B}$ denotes the $(r-k)$-logarithmic scale with center $(\Theta_k,...,\Theta_r)$ on $B$. We set
$$\mathcal{Y}_r(t,x):=(y_0(t,x),...,y_{k-1}(t,x),\hat{\mathcal{Y}}_{r-k,B}(t,\log(\vert{y_{k-1}(t,x)}\vert)))$$ 
for $(t,x) \in C$. With Definition \ref{2.1} one sees immediately that $\mathcal{Y}_r$ defines an $r$-logarithmic scale with center $\Theta:=(\Theta_0,...,\Theta_r)$ on $C$. Because $\Theta_j$ is $C$-nice for every $j \in \{0,...,r\}$ we see that $\mathcal{Y}_r$ is nice. In particular it is
$$\vert{\hat{\mathcal{Y}}_{r-k,B}(t,\log(\vert{y_{k-1}(t,x)}\vert))}\vert^{\otimes q}=\vert{\mathcal{Y}_r(t,x)}\vert^{\otimes q^*}$$
for every $q \in \mathbb{Q}^{r-k+1}$ and $(t,x) \in C$ where $q^*:=(0,...,0,q) \in \mathbb{Q}^{r+1}$. This observation gives the desired nice LA-preparing tuple for $\beta_j$ where $j \in \{1,...,m\}$.\hfill$\blacksquare$

\begin{rem}\label{2.41}
One may replace ''nice'' by ''pure'' and ''nicely'' by ''purely'' in Remark \ref{2.40}.
\end{rem}

\begin{defn}\label{2.42}
Let $r \in \mathbb{N}_0$. Let $m \in \mathbb{N}$. Let $g_1,...,g_m:C \to \mathbb{R}$ be functions. We call $g_1,...,g_m$ nicely $r$-log-analytically prepared \textbf{in a simultaneous way} if there is $\Theta:\pi(C) \to \mathbb{R}^{r+1}$ such that $g_1,...,g_m$ are nicely $r$-log-analytically prepared with center $\Theta$.
\end{defn}

\begin{prop}\label{2.43}
The following properties hold.
	\begin{itemize}
		\item[(1)]
		Let $r \in \mathbb{N}_0$. Let $g:C \to \mathbb{R}$ be nicely $r$-log-analytically prepared. Then there is $k \in \mathbb{N}$, a $C$-nice function $\eta:\pi(C) \to \mathbb{R}^k$, a globally subanalytic function $G:\mathbb{R}^k \times \mathbb{R}^{r+1} \to \mathbb{R}$, and a nice $r$-logarithmic scale $\mathcal{Y}_r$ on $C$ such that  
		$$g(t,x)=G(\eta(t),\mathcal{Y}_r(t,x))$$
		for all $(t,x) \in C$.
		
		\item[(2)]
		Let $r \in \mathbb{N}$. Let $h:C \to \mathbb{R}_{>0}$ be nicely $(r-1)$-log-analytically prepared. Then there is $k \in \mathbb{N}$, a $C$-nice function $\eta:\pi(C) \to \mathbb{R}^k$, a globally subanalytic function $H:\mathbb{R}^k \times \mathbb{R}^{r+1} \to \mathbb{R}$, a nice $(r-1)$-logarithmic scale $\mathcal{Y}_{r-1}:=(y_0,...,y_{r-1})$ on $C$ such that
		$$\log(h(t,x))=H(\eta(t),\mathcal{Y}_{r-1}(t,x),\log(\vert{y_{r-1}(t,x)}\vert))$$
		for all $(t,x) \in C$.
	\end{itemize}
\end{prop}

{\bf Proof:}
$(1)$: Let
$$(r,\mathcal{Y}_r,a,q,s,v,b,P)$$ 
be a nice LA-preparing tuple for $g$. Take $k:=s+1$,
$$\eta=(\eta_1,...,\eta_k):\pi(C) \to \mathbb{R}^k, t \mapsto (a(t),b_1(t),...,b_s(t)).$$
Then $\eta$ is $C$-nice. Let $z:=(z_0,...,z_s)$ range over $\mathbb{R}^{s+1}$ and $w:=(w_0,...,w_r)$ range over $\mathbb{R}^{r+1}$. Set
$$\alpha_0:\mathbb{R}^k \times \mathbb{R}^{r+1} \to \mathbb{R}, (z,w) \mapsto z_0 \prod_{j=0}^r \vert{w_j}\vert^{q_j}.$$ 
For $i \in \{1,...,s\}$ let
$$\alpha_i:\mathbb{R}^k \times \mathbb{R}^{r+1} \to \mathbb{R}, (z,w) \mapsto z_i \prod_{j=0}^r \vert{w_j}\vert^{p_{ij}}.$$
Set
$$G:\mathbb{R}^k \times \mathbb{R}^{r+1} \to \mathbb{R}, (z_0,...,z_s,w_0,....,w_r) \mapsto $$
$$\left\{\begin{array}{ll} \alpha_0(z,w)v(\alpha_1(z,w),...,\alpha_s(z,w)), & \vert{\alpha_i(z,w)}\vert \leq 1 \textnormal{ for all } i \in \{1,...,s\}, \\
0, & \textnormal{otherwise}. \end{array}\right.$$
Then $G$ is globally subanalytic and for every $(t,x) \in C$ we have
$$g(t,x)=G(\eta(t),\mathcal{Y}_r(t,x)).$$

$(2)$: Let
$$(r-1,\mathcal{Y}_{r-1},a,q,s,v,b,P)$$ 
be a nice LA-preparing tuple for $h$ where $\mathcal{Y}_{r-1}:=(y_0,...,y_{r-1})$ is a nice $(r-1)$-logarithmic scale with center $(\Theta_0,...,\Theta_{r-1})$ on $C$. Then $a>0$ on $C$. Take $k:=s+r+1$,
$$\eta=(\eta_1,...,\eta_k):\pi(C) \to \mathbb{R}^k, t \mapsto$$
$$(\log(a(t)),b_1(t),...,b_s(t),\Theta_1(t),...,\Theta_{r-1}(t),0).$$

Then $\eta$ is $C$-nice by Remark \ref{2.35}. Let $z:=(z_0,...,z_{s+r})$ range over $\mathbb{R}^{s+r+1}$ and $w:=(w_0,...,w_r)$ over $\mathbb{R}^{r+1}$. Set 
$$\beta:\mathbb{R}^k \times \mathbb{R}^{r+1} \to \mathbb{R}, (z,w) \mapsto z_0 + \sum_{j=0}^{r-1}q_j(w_{j+1}+z_{s+j+1}).$$
For $i \in \{1,...,s\}$ let
$$\alpha_i:\mathbb{R}^k \times \mathbb{R}^{r+1} \to \mathbb{R}, (z,w) \mapsto z_i \prod_{j=0}^{r-1} \vert{w_j}\vert^{p_{ij}}.$$
Set
$$H:\mathbb{R}^k \times \mathbb{R}^{r+1} \to \mathbb{R}, (z_0,....,z_{s+r},w_0,....,w_r) \mapsto$$
$$\left\{\begin{array}{ll} \beta(z,w) + \log(v(\alpha_1(z,w),...,\alpha_s(z,w))), & \vert{\alpha_i(z,w)}\vert \leq 1 \textnormal{ for all } i \in \{1,...,s\}, \\
0, & \textnormal{otherwise}. \end{array}\right.$$
Then $H$ is globally subanalytic since $\log(v)$ is globally subanalytic. For every $(t,x) \in C$ we have
$$\log(h(t,x))=H(\eta(t),\mathcal{Y}_{r-1}(t,x),\log(\vert{y_{r-1}(t,x)}\vert)).$$
\hfill$\blacksquare$

\begin{cor}\label{2.44}
Let $r,k,l \in \mathbb{N}$. Let $g_1,...,g_k:C \to \mathbb{R}$ and $h_1,...,h_l:C \to \mathbb{R}_{>0}$ be nicely $(r-1)$-log-analytically prepared in a simultaneous way. Let $F:\mathbb{R}^{k+l} \to \mathbb{R}$ be globally subanalytic. Then there are $m \in \mathbb{N}$, a $C$-nice $\eta:\pi(C) \to \mathbb{R}^m$, a globally subanalytic $J:\mathbb{R}^m \times \mathbb{R}^{r+1} \to \mathbb{R}$, and a nice $(r-1)$-logarithmic scale $\mathcal{Y}:=(y_0,...,y_{r-1})$ on $C$ such that
for all $(t,x) \in C$
$$F(g_1(t,x),...,g_k(t,x),\log(h_1(t,x)),...,\log(h_l(t,x)))$$
$$=J(\eta(t),\mathcal{Y}(t,x),\log(\vert{y_{r-1}(t,x)}\vert)).$$
\end{cor}

\begin{rem}\label{2.45}
One may replace ''nicely prepared'' by ''purely prepared'', ''$C$-nice'' by ''log-analytic'' and ''nice $r$-logarithmic scale'' by '' pure $r$-logarithmic scale'' in Proposition \ref{2.43} and Corollary \ref{2.44}.
\end{rem}

Here is the promised class of log-analytic functions in $x$ which is closed under log-analytic preparation.

\begin{defn}\label{2.46}
Let $X \subset \mathbb{R}^n \times \mathbb{R}$ be definable and let $f:X \to \mathbb{R}$ be a function.
\begin{itemize}
	\item[(a)] Let $r \in \mathbb{N}_0$. By induction on $r$ we define that $f$ is \textbf{nicely log-analytic in $x$ of order at most $r$}.
	
	\textbf{Base case}: The function $f$ is log-analytic in $x$ of order at most $0$ if the following holds: There is a definable cell decomposition $\mathcal{C}$ of $X$ such that for every $C \in \mathcal{C}$ there is $m \in \mathbb{N}$, a globally subanalytic $F:\mathbb{R}^m \times \mathbb{R} \to \mathbb{R}$ and a $C$-nice $g:\pi(C) \to \mathbb{R}^m$ such that $f(t,x)=F(g(t),x)$ for every $(t,x) \in C$.
	
	\textbf{Inductive step}: The function $f$ is nicely log-analytic in $x$ of order at most $r$ if the following holds: There is a definable cell decomposition $\mathcal{C}$ of $X$ such that for $C \in \mathcal{C}$ there are $k,l \in \mathbb{N}_0$, nicely log-analytic functions $g_1,...,g_k:C \to \mathbb{R}, h_1,...,h_l:C \to \mathbb{R}_{>0}$ in $x$ of order at most $r-1$, and a globally subanalytic $F:\mathbb{R}^{k+l} \to \mathbb{R}$ such that
	$$f|_C=F(g_1,...,g_k,\log(h_1),...,\log(h_l)).$$
	\item[(b)] Let $r \in \mathbb{N}_0$. The function $f$ is \textbf{nicely log-analytic in $x$ of order $r$} if $f$ is nicely log-analytic in $x$ of order at most $r$ but not nicely log-analytic in $x$ of order at most $r-1$.
	\item[(c)] The function $f$ is \textbf{nicely log-analytic in $x$} if there is $r \in \mathbb{N}_0$ such that $f$ is nicely log-analytic in $x$ of order $r$.
\end{itemize}
\end{defn}

\begin{rem}\label{2.47}
Let $r \in \mathbb{N}_0$. Then the following properties hold.
\begin{itemize}
\item [(1)] Let $f:C \to \mathbb{R}$ be a function. If $f$ is purely log-analytic in $x$ of order at most $r$ then it is nicely log-analytic in $x$ of order at most $r$.
\item [(2)] Let $m \in \mathbb{N}$. Let $f_1,...,f_m:C \to \mathbb{R}$ be nicely log-analytic in $x$ of order at most $r$ and let $F:\mathbb{R}^m \to \mathbb{R}$ be globally subanalytic. Then $F(f_1,...,f_m)$ is nicely log-analytic in $x$ of order at most $r$.
\item [(3)] Let $m \in \mathbb{N}$. Let $f_1,...,f_m:C \to \mathbb{R}$ be nicely log-analytic in $x$ and let $F:\mathbb{R}^m \to \mathbb{R}$ be log-analytic. Then $F(f_1,...,f_m)$ is nicely log-analytic in $x$.
\end{itemize}
\end{rem}

\begin{rem}\label{2.48}
A nice logarithmic scale $\mathcal{Y}:=(y_0,...,y_r)$ is nicely log-analytic. In particular $y_j$ is nicely log-analytic in $x$ of order at most $j$ for $j \in \{0,...,r\}$.
\end{rem}

{\bf Proof:}
Let $\mathcal{Y}:=(y_0,...,y_r)$. Let $(\Theta_0,...,\Theta_r)$ be the $C$-nice center of $\mathcal{Y}$. We show by induction on $j \in \{0,...,r\}$ that $y_j$ is nicely log-analytic in $x$ of order at most $j$.

$j=0$: We have $y_0(t,x)=x-\Theta_0(t)$ for every $(t,x) \in C$. That $y_0$ is nicely log-analytic in $x$ of order $0$ is clear with (2) in Remark \ref{2.47}.

$j-1 \to j$: It holds $y_j(t,x)=\log(\vert{y_{j-1}(t,x)}\vert)-\Theta_j(t)$ for every $(t,x) \in C$. That $y_j$ is nicely log-analytic in $x$ of order at most $j$ is clear with (2) in Remark \ref{2.47}. \hfill$\blacksquare$

\begin{prop}\label{2.49}
Let $k \in \mathbb{N}_0$ and $r \in \mathbb{N}$. Let $F:\mathbb{R}^{k+r+1} \to \mathbb{R}$ be globally subanalytic and $\eta:\pi(C) \to \mathbb{R}^k$ be $C$-nice. Let $\mathcal{Y}:=(y_0,...,y_r)$ be a nice $r$-logarithmic scale with center $(\Theta_0,...,\Theta_r)$ on $C$. Let
$$f:C \to \mathbb{R}, (t,x) \mapsto F(\eta(t),\mathcal{Y}(t,x)).$$
Then the following properties hold.
\begin{itemize}
\item [(1)] $f$ is nicely log-analytic in $x$ of order at most $r$.
\item [(2)] Assume that there is a $C$-nice $\xi:\pi(C) \to \mathbb{R}$ such that $y_0 \sim_C \xi$. Then $f$ is nicely log-analytic in $x$ of order at most $r-1$.
\end{itemize}
\end{prop}

{\bf Proof:}
(1): By Remark \ref{2.48} $y_l$ is nicely log-analytic in $x$ of order at most $l$ for every $l \in \{0,...,r\}$. We are done with (2) in Remark \ref{2.47}.

(2): By (1) in Remark \ref{2.11} there is $\delta>1$ such that $1/\delta<\frac{y_0}{\xi}<\delta$ on $C$. Set
$$
\log^*:\mathbb{R} \to \mathbb{R}, y \mapsto \left\{\begin{array}{ll} \log(y), & y \in [1/\delta,\delta], \\
0, & y \not\in [1/\delta,\delta]. \end{array}\right.
$$
Then $\log^*$ is globally subanalytic. We have $y_1=y_1^*$ on $C$ where
$$y_1^*:= \log^*\bigl(\frac{y_0}{\xi}\bigl)+ \log(\vert{\xi}\vert) - \Theta_1.$$
Then $y_1^*$ is nicely log-analytic in $x$ of order $0$, because $\log(\vert{\xi}\vert)$ is $C$-nice by Remark \ref{2.35}. In particular 
$$f(t,x)=F(\eta(t),y_0(t,x),y^*_1(t,x),...,y^*_r(t,x))$$
for every $(t,x) \in C$ where inductively for $l \in \{2,...,r-1\}$
$$y_l^*:= \log(\vert{y^*_{l-1}}\vert) - \Theta_l.$$
Note that $y_l^*$ is nicely log-analytic in $x$ of order at most $l-1$ for every $l \in \{2,...,r\}$ and therefore $F(\eta,y_0,y^*_1,...,y^*_r)$ is nicely log-analytic in $x$ of order at most $r-1$ by (2) in Remark \ref{2.47}. \hfill$\blacksquare$

\begin{rem}\label{2.50}
Let $r \in \mathbb{N}_0$. Let $f:C \to \mathbb{R}$ be a function. If $f$ is nicely $r$-log-analytically prepared in $x$ then $f$ is nicely log-analytic in $x$ of order at most $r$.
\end{rem}

{\bf Proof:}
By (1) in Proposition \ref{2.43} there is $k \in \mathbb{N}$, a $C$-nice function $\eta:\pi(C) \to \mathbb{R}^k$, a globally subanalytic function $G:\mathbb{R}^k \times \mathbb{R}^{r+1} \to \mathbb{R}$, and a nice $r$-logarithmic scale $\mathcal{Y}$ on $C$ such that 
$$g(t,x)=G(\eta(t),\mathcal{Y}(t,x))$$
for all $(t,x) \in C$. By (1) in Proposition \ref{2.49} $g$ is nicely log-analytic in $x$ of order at most $r$. \hfill$\blacksquare$

\begin{prop}\label{2.51}
Let $k \in \mathbb{N}_0$ and $r \in \mathbb{N}$. Let $F:\mathbb{R}^{k+r} \to \mathbb{R}$ be globally subanalytic and $\eta:\pi(C) \to \mathbb{R}^k$ be $C$-nice. Let $\mathcal{Y}:=(y_0,...,y_r)$ be a nice $r$-logarithmic scale with center $(\Theta_0,...,\Theta_r)$ on $C$. Let
$$f:C \to \mathbb{R}, (t,x) \mapsto F(\eta(t),y_1(t,x),...,y_r(t,x)).$$
Then the following holds.
\begin{itemize}
	\item [(1)]
	There is a nicely log-analytic function $\kappa:C^{\log(\vert{y_0}\vert)} \to \mathbb{R}$ in $x$ of order at most $r-1$ such that $f(t,x)=\kappa(t,\log(\vert{y_0(t,x)}\vert))$ for every $(t,x) \in C$.
	\item [(2)]
	Assume $r \geq 2$. Let $l \in \{1,...,r-1\}$. Let $\xi:\pi(C) \to \mathbb{R}$ be $C$-nice such that $y_l \sim_C \xi$. Then there is a nicely log-analytic function $\lambda:C^{\log(\vert{y_0}\vert)} \to \mathbb{R}$ in $x$ of order at most $r-2$ such that $f(t,x)=\lambda(t,\log(\vert{y_0(t,x)}\vert))$ for every $(t,x) \in C$. 
\end{itemize}
\end{prop}

{\bf Proof:}
We set $B:=C^{\log(\vert{y_0}\vert)}$. Let $\mu_0:B \to \mathbb{R}, (t,x) \mapsto x-\Theta_1(t)$ and inductively for $j \in \{1,...,r-1\}$ let $\mu_j:B \to \mathbb{R}, (t,x) \mapsto \log(\vert{\mu_{j-1}(t,x)}\vert)-\Theta_{j+1}(t)$. Note that 
$$y_j(t,x)=\mu_{j-1}(t,\log(\vert{y_0(t,x)}\vert))$$
for every $(t,x) \in C$ and $j \in \{1,...,r\}$. With Remark \ref{2.5} we obtain that $\mathcal{Y}_{r-1,B}:=(\mu_0,....,\mu_{r-1})$ is a nice $(r-1)$-logarithmic scale with center $(\Theta_1,...,\Theta_r)$ on $B$. It holds $$f(t,x)=F(\eta(t),\mathcal{Y}_{r-1,B}(t,\log(\vert{y_0(t,x)}\vert)))$$ 
for every $(t,x) \in C$.

(1): We set 
$$\kappa:B \to \mathbb{R}, (t,x) \mapsto F(\eta(t),\mathcal{Y}_{r-1,B}(t,x)).$$
By (1) in Proposition \ref{2.49} $\kappa$ is a nicely log-analytic function in $x$ of order at most $r-1$ and we obtain
$$f(t,x)=\kappa(t,\log(\vert{y_0(t,x)}\vert))$$
for all $(t,x) \in C$.

(2): Since $y_l \sim_C \xi$ we obtain $\mu_{l-1} \sim_B \xi$ since $\xi$ depends only on $t$ and $\pi(B)=\pi(C)$. Thus by (1) in Remark \ref{2.11} there is $\delta>1$ such that $1/\delta <\tfrac{\mu_{l-1}}{\xi}<\delta$ on $B$. Consider
$$\log^*:\mathbb{R} \to \mathbb{R}, y \mapsto \left\{\begin{array}{ll} \log(y), & y \in [1/\delta,\delta], \\
0, & y \not\in [1/\delta,\delta]. \end{array}\right.$$
Then $\log^*$ is globally subanalytic. Set
$$\mu_l^*:=\log^*\Bigl(\frac{\mu_{l-1}}{\xi}\Bigl)+\log(\vert{\xi}\vert)-\Theta_{l+1}$$
and inductively for $j \in \{l+1,...,r-1\}$ 
$$\mu_j^*:= \log(\vert{\mu_{j-1}^*}\vert)-\Theta_{j+1}.$$
Then by Remark \ref{2.48} and (2) in \ref{2.47} $\mu_l^*$ is a nicely log-analytic function in $x$ of order at most $l-1$, because $\log(\vert{\xi}\vert)$ is $C$-nice by Remark \ref{2.35}. We see similarly as in the proof of Remark \ref{2.48} that $\mu_j^*$ is nicely log-analytic in $x$ of order at most $j-1$ for every $j \in \{l+1,...,r-1\}$. We set
$$\lambda:B \to \mathbb{R}, (t,x) \mapsto F(\eta(t),\mu_0(t,x),...,\mu_{l-1}(t,x),\mu_l^*(t,x),...,\mu_{r-1}^*(t,x)).$$

By (2) in Remark \ref{2.47} $\lambda$ is nicely log-analytic in $x$ of order at most $r-2$.
Because $\mu_j=\mu^*_j$ for every $j \in \{l,...,r-1\}$ on $B$ we obtain
$$f(t,x)=\lambda(t,\log(\vert{y_0(t,x)}\vert))$$
for every $(t,x) \in C$. \hfill$\blacksquare$

An immediate consequence from the globally subanalytic preparation theorem is the following.

\begin{prop}\label{2.52}
Let $X \subset \mathbb{R}^n \times \mathbb{R}$ be definable. Let $m \in \mathbb{N}$. Let $f_1,...,f_m: X \to \mathbb{R}$ be nicely log-analytic functions in $x$ of order $0$. Then there is a definable cell decomposition $\mathcal{C}$ of $X_{\neq 0}$ such that $f_1|_C,...,f_m|_C$ are nicely $0$-log-analytically prepared in a simultaneous way for every $C \in \mathcal{C}$.
\end{prop}
%auf neq 0 kann nicht verzichtet werden!

{\bf Proof:} It is enough to consider the following situation: Let $g_1,...,g_k:\pi(X) \to \mathbb{R}$ be $X$-nice and $F_j:\mathbb{R}^{k+1} \to \mathbb{R}$ be globally subanalytic  such that
$$f_j(t,x) = F_j(g_1(t),...,g_k(t),x)$$
for every $(t,x) \in X$ and $j \in \{1,...,m\}$. Let $g:=(g_1,...,g_k)$. Let $z:=(z_1,...,z_k)$ range over $\mathbb{R}^k$. Let $\pi^*:\mathbb{R}^{k+1}\to \mathbb{R}^k, (z,x) \mapsto z$. With Fact 2.15 we find a globally subanalytic cell decomposition $\mathcal{D}$ of $\mathbb{R}^k \times \mathbb{R}_{\neq 0}$ such that $F_1|_D,...,F_m|_D$ are globally subanalytically prepared in $x$ in a simultaneous way for every $D \in \mathcal{D}$. There is a definable cell decomposition $\mathcal{C}$ of $X_{\neq 0}$ such that for every $C \in \mathcal{C}$ there is $D_C \in \mathcal{D}$ such that $(g(t),x) \in D_C$ for every $(t,x) \in C$. Fix a $C \in \mathcal{C}$, the globally subanalytic center $\vartheta$, and for $j \in \{1,...,m\}$ a preparing tuple $(0,y_0,a_j,q_j,s,v_j,b,P)$ for $F_j|_{D_C}$ where $y_0:=y-\vartheta(z)$ on $D_C$, $b:=(b_1,...,b_s)$, $P:=(p_1,...,p_s)^t$, and $a_1,...,a_m,b_1,...,b_s$ are globally subanalytic on $\pi^*(D_C)$. We have for every $j \in \{1,...,m\}$ and every $(z,x) \in D_C$
$$F_j(z,x)=a_j(z) \vert{x-\vartheta(z)}\vert^{q_j} v_j(b_1(z)\vert{x-\vartheta(z)}\vert^{p_1},...,b_s(z)\vert{x-\vartheta(z)}\vert^{p_s}).$$ 
Let $h:C \to \mathbb{R}, (t,x) \mapsto x-\vartheta(g(t))$. Then it is immediately seen with Definition \ref{2.14} that $h$ is a $0$-logarithmic scale. We obtain
$$f_j(t,x)=a_j(g(t)) \vert{h(t,x)}\vert^{q_j}v_j(b_1(g(t))\vert{h(t,x)}\vert^{p_1},...,b_s(g(t))\vert{h(t,x)}\vert^{p_s})$$
for every $(t,x) \in C$ and are done, because $a_1(g),...,a_m(g),\vartheta(g),b_1(g),...,b_s(g):\pi(C) \to \mathbb{R}$ are $C$-nice by Remark \ref{2.35}. \hfill$\blacksquare$

\begin{rem}\label{2.53}
One may replace ''nicely log-analytic'' by ''purely log-analytic'' and ''nicely $0$-log-analytically prepared'' by ''purely $0$-log-analytically prepared'' in Proposition \ref{2.52}.
\end{rem}

\begin{prop}\label{2.54}
Let $m \in \mathbb{N}$ and $r \in \mathbb{N}_0$. Let $X \subset \mathbb{R}^n \times \mathbb{R}$ be definable. Let $f_1,....,f_m:X \to \mathbb{R}$ be nicely log-analytic functions in $x$ of order at most $r$. Then there is a definable cell decomposition $\mathcal{C}$ of $X_{\neq 0}$ such that $f_1|_C,... ,f_m|_C$ are nicely $r$-log-analytically prepared in a simultaneous way for every $C \in \mathcal{C}$.
\end{prop}
%auf neq 0 kann nicht verzichtet werden, da r global ist!

{\bf Proof:}
We do an induction on $r$. 

$r=0$: Then $f_1,...,f_m$ are nicely log-analytic in $x$ of order $0$ and we are done with Proposition \ref{2.52}. 

$<r \to r$: It is enough to consider the following situation: Assume that there are $l,l' \in \mathbb{N}_0$, nicely log-analytic functions $g_1,...,g_l:X \to \mathbb{R}$ and $h_1,...,h_{l'}:X \to \mathbb{R}_{>0}$ in $x$ of order at most $r-1$ and for every $j \in \{1,...,m\}$ there is a globally subanalytic function $F_j:\mathbb{R}^{l+l'} \to \mathbb{R}$ such that
$$f_j = F_j(g_1,...,g_l,\log(h_1),..,\log(h_{l'})).$$
Applying the inductive hypothesis to $g_1,...,g_l,h_1,...,h_{l'}$ and Corollary \ref{2.44} we find a definable cell decomposition $\mathcal{U}$ of $X_{\neq 0}$ such that for every $U \in \mathcal{U}$ there is a nice $(r-1)$-logarithmic scale $\mathcal{Y}_{r-1}:=(y_0,...,y_{r-1})$ on $U$, a $k \in \mathbb{N}$, a $U$-nice function $\eta:\pi(U) \to \mathbb{R}^k$, and for every $j \in \{1,...,m\}$ there is a globally subanalytic function $H_j:\mathbb{R}^k \times \mathbb{R}^{r+1} \to \mathbb{R}$ such that
$$f_j(t,x)=H_j(\eta(t),\mathcal{Y}_{r-1}(t,x),\log(\vert{y_{r-1}(t,x)}\vert))$$
for all $(t,x) \in U$. Fix $U \in \mathcal{U}$ and for this $U$ a corresponding $\mathcal{Y}_{r-1}$, $\eta$, and $H_j$ for $j \in \{1,...,m\}$. By further decomposing $U$ if necessary we may assume that either $\vert{y_{r-1}}\vert=1$ or $\vert{y_{r-1}}\vert>1$ or $\vert{y_{r-1}}\vert<1$ on $U$. Assume the former. Then by (1) in Proposition \ref{2.49} $f|_C$ is nicely log-analytic in $x$ of order at most $(r-1)$ and we are done with the inductive hypothesis. Assume $\vert{y_{r-1}}\vert>1$ or $\vert{y_{r-1}}\vert<1$ on $U$. Then $\mathcal{Y}:=(y_0,...,y_r)$ is an $r$-logarithmic scale on $U$ where $y_r:=\log(\vert{y_{r-1}}\vert)$.

Let $(z,w):=(z_1,...,z_k,w_0,...,w_r)$ range over $\mathbb{R}^k \times \mathbb{R}^{r+1}$. Set $w':=(w_1,...,w_r)$. Let $\pi^*:\mathbb{R}^k \times \mathbb{R}^{r+1} \to \mathbb{R}^{r+k}$ be the projection on $(z_1,...,z_k,w_1,...,w_r)$. With Fact 2.15 we find a globally subanalyic cell decomposition $\mathcal{D}$ of $\mathbb{R}^k \times \mathbb{R}^r \times \mathbb{R}_{\neq 0}$ such that $H_1|_D,...,H_m|_D$ are globally subanalytically prepared in $w_0$ in a simultaneous way. There is a definable cell decomposition $\mathcal{A}$ of $U$ such that for every $A \in \mathcal{A}$ there is $D_A \in \mathcal{D}$ such that $(\eta(t),\mathcal{Y}_r(t,x)) \in D_A$ for every $(t,x) \in A$. Fix $A \in \mathcal{A}$, the globally subanalytic center $\vartheta$, and for $j \in \{1,...,m\}$ an LA-preparing tuple 
$(0,y,a_j,q_j,s,v_j,b,P)$ for $F_j|_{D_A}$ where $y:=w_0-\vartheta(w',z)$ on $D_A$, $b:=(b_1,...,b_s)$, $P:=(p_1,...,p_s)^t$, and $a_1,...,a_m,b_1,...,b_s$
are globally subanalytic on $\pi^*(D_A)$. We have for every $j \in \{1,...,m\}$ 
$$H_j|_{D_A}=a_j(z,w')\vert{y(z,w)}\vert^{q_j} v_j(b_1(z,w')\vert{y(z,w)}\vert^{p_1},...,b_s(z,w')\vert{y(z,w)}\vert^{p_s}).$$
%s unabhängig von j, ergänze die Reihe!

From the inductive hypothesis we will derive the following two claims.

{\bf Claim 1}\\
Let $d \in \mathbb{N}$ and $p \in \{1,...,r\}$. Let $\alpha_1,...,\alpha_d:A^{\log(\vert{y_0}\vert)} \to \mathbb{R}$ be nicely log-analytic in $x$ of order at most $p-1$. For $j \in \{1,...,d\}$ consider
$$\beta_j:A \to \mathbb{R}, (t,x) \mapsto \alpha_j(t,\log(\vert{y_0(t,x)}\vert)).$$
Then there is a definable cell decomposition $\mathcal{Q}$ of $A$ such that for every $Q \in \mathcal{Q}$ the following holds. 
\begin{itemize}
	\item [(1)] There are $\hat{\Theta}_1,...,\hat{\Theta}_p:\pi(Q) \to \mathbb{R}$ such that $\beta_1|_Q,...,\beta_d|_Q$ are nicely $p$-log-analytically prepared with center $\hat{\Theta}:=(\Theta_0|_{\pi(Q)},\hat{\Theta}_1,...,\hat{\Theta}_p)$.
	\item [(2)] Additionally for every $j \in \{1,...,d\}$ there is for $\beta_j|_Q$ a nice LA-preparing tuple $(p,\hat{\mathcal{Y}}_p,\hat{a}_j,\hat{q}_j,\hat{s},\hat{v}_j,\hat{b},\hat{P})$ such that $\hat{q}_j \in \{0\} \times \mathbb{Q}^p$ and $\hat{P} \in M_1(\hat{s} \times (p+1),\mathbb{Q})$ where $\hat{\mathcal{Y}}_p$ denotes the $p$-logarithmic scale with center $\hat{\Theta}$ on $Q$.
\end{itemize}

\textbf{Proof of Claim 1}

Set $B:=A^{\log(\vert{y_0}\vert)}$. Applying the inductive hypothesis to $\alpha_1,...,\alpha_d$ we obtain a definable cell decomposition $\mathcal{S}$ of $B$ such that $\alpha_1 |_S,...,\alpha_d |_S$ are nicely $(p-1)$-log-analytically prepared in a simultaneous way for every $S \in \mathcal{S}$. Consider the definable set
$$S^{\log{(\vert{y_0}\vert)}^*}:=\{(t,x) \in A \mid (t,\log(\vert{y_0(t,x)}\vert)) \in S\}$$
for $S \in \mathcal{S}$. It holds $\bigcup_{S \in \mathcal{S}} S^{\log{(\vert{y_0}\vert)}^*}=A$. Fix $S \in \mathcal{S}$ and the center $(\hat{\Theta}_1,...,\hat{\Theta}_p)$ of $\alpha_1|_S,...,\alpha_d|_S$. Let $T:=S^{{\log(\vert{y_0}\vert)}^*}$. Note that $\beta_1|_T,...,\beta_d|_T$ are nicely $p$-log-analytically prepared with center $(\Theta_0|_{\pi(T)},\hat{\Theta}_1,...,\hat{\Theta}_p)$ and for $j \in \{1,...,d\}$ there is a nice LA-preparing tuple $(p,\hat{\mathcal{Y}}_{p,T},\hat{a}_j,\hat{q}_j,\hat{s},\hat{v}_j,\hat{b},\hat{P})$ for $\beta_j|_T$ such that $\hat{q}_j \in \{0\} \times \mathbb{Q}^p$ and $\hat{P} \in M_1(\hat{s} \times (p+1), \mathbb{Q})$ by Remark \ref{2.40} where $\hat{\mathcal{Y}}_{p,T}$ denotes the nice $p$-logarithmic scale with center $(\Theta_0|_{\pi(T)},\hat{\Theta}_1,...,\hat{\Theta}_p)$ on $T$. With the cell decomposition theorem applied to $T$ we are done. \hfill $\blacksquare_{\textnormal{Claim }1}$

{\bf Claim 2}\\
Let $\gamma:\mathbb{R}^k \times \mathbb{R}^r \to \mathbb{R}$ be globally subanalytic such that 
$$y_0 \sim_A \gamma(\eta,y_1,...,y_r).$$ 
Then $f_j|_A$ is nicely log-analytic in $x$ of order at most $r-1$ for every $j \in \{1,...,m\}$.

{\bf Proof of Claim 2} By (1) in Proposition \ref{2.51} there is a nicely log-analytic function $\kappa:A^{\log(\vert{y_0}\vert)} \to \mathbb{R}$ in $x$ of order at most $r-1$ such that $$\gamma(\eta,y_1,...,y_r)=\kappa(t,\log(\vert{y_0}\vert))$$
on $A$. So by Claim $1$ there is a definable cell decomposition $\mathcal{N}$ of $A$ such that for every $N \in \mathcal{N}$ the function $\gamma(\eta,y_1,...,y_r)|_N$ is nicely $r$-log-analytically prepared and there is a nice LA-preparing tuple $(r,\tilde{\mathcal{Y}}_r,\Psi,\tilde{q},\tilde{s},\tilde{v},\tilde{b},\tilde{P})$ for $\gamma(\eta,y_1,...,y_r)|_N$ such that $\tilde{q} \in \{0\} \times \mathbb{Q}^r$ and $\tilde{P} \in M_1(\tilde{s} \times (r+1), \mathbb{Q})$ where $\tilde{\mathcal{Y}}_r$ has a center whose first component is $\Theta_0|_{\pi(N)}$. Fix $N \in \mathcal{N}$. 

By Remark \ref{2.18} it is enough to consider the following property $(*)_p$ for $p \in \{0,...,r\}$ on $N$: There is a nice $p$-logarithmic scale $\tilde{\mathcal{Y}}_p:=(y_0,\tilde{y}_1,...,\tilde{y}_p)$ with center $(\Theta_0|_{\pi(N)},\tilde{\Theta}_1,...,\tilde{\Theta}_p)$, and an $N$-nice $\Psi:\pi(N) \to \mathbb{R}$ such that
$$y_0 \sim_N \Psi \prod_{l=1}^p \vert{\tilde{y}_l}\vert^{\tilde{q}_l}$$
where $\tilde{q}_l \in \mathbb{Q}$ for every $l \in \{1,...,p\}$.

If $p=0$ then $y_0 \sim_N \Psi$. We are done with (2) in Proposition \ref{2.49} applied to $f_j|_N=H_j(\eta|_N,\mathcal{Y}|_N)$ for every $j \in \{1,...,m\}$. Assume $p>0$. By a suitable induction on $p$ it is enough to establish the following Subclaim.

{\bf Subclaim}\\
There is a decomposition $\mathcal{K}$ of $N$ into finitely many definable sets such that the following holds for every $K \in \mathcal{K}$: The function $f_j|_K$ is nicely log-analytic in $x$ of order at most $r-1$ for every $j \in \{1,...,m\}$ or $(*)_{p-1}$ holds on $K$. 

{\bf Proof of the Subclaim:} For $M>1$ let $N_{>M}:=N_{>M}(\tilde{\mathcal{Y}}_p)$ and for $i \in \{1,...,p\}$ let $N_{i,M}:=N_{i,M}(\tilde{\mathcal{Y}}_p)$. By Remark \ref{2.38} there is $M>1$ and an $N$-nice $\xi:\pi(N) \to \mathbb{R}$ such that $y_0 \sim_{N_{>M}} \xi$. Fix such an $M$. Since 
$$N=N_{>M} \cup N_{1,M} \cup ... \cup N_{p,M}$$   
it suffices to establish the Subclaim for $N_{>M}$ and $N_{i,M}$ instead for $N$ where $i \in \{1,...,p\}$.

$B:=N_{>M}$: By (2) in Remark \ref{2.34} $\xi|_B$ is $B$-nice. Therefore with (2) in Proposition \ref{2.49} applied to $f_j|_B=H_j(\eta|_B,\mathcal{Y}_r|_B)$
we obtain that $f_j|_B$ is a nicely log-analytic function in $x$ of order at most $r-1$.

$B_i:=N_{i,M}$ for $i \in \{1,...,p\}$: Set $\tilde{y}_0:=y_0$. It holds $\vert{\tilde{y}_i}\vert<M$ on $B_i$. So we obtain
$$\frac{1}{\delta} < \frac{\vert{\tilde{y}_{i-1}}\vert}{\exp(\tilde{\Theta}_i)}<\delta$$
for $\delta:=e^M$ on $B_i$ which gives $\tilde{y}_{i-1} \sim_{B_i} \exp(\tilde{\Theta}_i)|_{\pi(B_i)}$.

Assume $i=1$. With Remark \ref{2.33} we see that $\exp(\tilde{\Theta}_1)|_{\pi(B_1)}$ is a $B_1$-nice $B_1$-heir. Again with (2) in Proposition \ref{2.49} applied to $f_j|_{B_1}=H_j(\eta|_{B_1},\mathcal{Y}|_{B_1})$ (we have $y_0 \sim_{B_1} \xi$ with $\xi=\exp(\tilde{\Theta}_1)|_{\pi(B_1)}$) we see that $f_j|_{B_1}$ is nicely log-analytic in $x$ of order at most $r-1$ for every $j \in \{1,...,m\}$.

Assume $i>1$. Then $p>1$. Let 
$$\Xi:N \to \mathbb{R}, (t,x) \mapsto \Psi(t) \prod_{l=1}^p\vert{\tilde{y}_l(t,x)}\vert^{\tilde{q}_l}.$$
By (the proof of) (1) in Proposition \ref{2.43} there is a globally subanalytic function $G:\mathbb{R} \times \mathbb{R}^p \to \mathbb{R}$ such that $\Xi=G(\Psi,\tilde{y}_1,...,\tilde{y}_p)$ on $B_i$. By (2) in Proposition \ref{2.51} there is a nicely log-analytic function $\lambda:B_i^{\log{(\vert{y_0}\vert)}} \to \mathbb{R}$ in $x$ of order at most $p-2$ such that 
$$\Xi(t,x)=\lambda(t,\log(\vert{y_0(t,x)}\vert))$$
for every $(t,x) \in B_i$. Note that $0 \leq p-2 < r$. With Claim 1 applied to $\Xi$ we find a definable cell decomposition $\mathcal{K}$ of $B_i$ such that $\Xi|_K$ is nicely $(p-1)$-log-analytically prepared with nice LA-preparing tuple $(p-1,\mathcal{Y}_{p-1},\bar{a}_j,\bar{q}_j,\bar{s},\bar{v}_j,\bar{b},\bar{P})$ where $\bar{q}_j \in \{0\} \times \mathbb{Q}^{p-1}$, $\bar{P} \in M_1(\bar{s} \times p,\mathbb{Q})$, and $\mathcal{Y}_{p-1}$ is a nice $(p-1)$-logarithmic scale with center $\bar{\Theta}:=(\Theta_0|_{\pi(K)},\bar{\Theta}_1,...,\bar{\Theta}_{p-1})$. With Remark \ref{2.18} we obtain property $(*)_{p-1}$ applied to every $K \in \mathcal{K}$. \hfill$\blacksquare_{\textnormal{Subclaim}}$

\hfill$\blacksquare_{\textnormal{Claim }2}$

\textbf{Case 1:} $\vartheta=0$. Let $a_{m+j}:=b_j$ for $j \in \{1,...,s\}$. By (1) in Proposition \ref{2.51} there are nicely log-analytic functions $\alpha_1,...,\alpha_{m+s}:A^{\log(\vert{y_0}\vert)} \to \mathbb{R}$  in $x$ of order at most $r-1$ such that for every $(t,x) \in A$ we have
$$a_j(\eta(t),y_1(t,x),...,y_r(t,x))=\alpha_j(t,\log(\vert{y_0(t,x)}\vert))$$
for every $j \in \{1,...,m+s\}$. With Claim 1 applied to 
$$\beta_j:A \to \mathbb{R}, (t,x) \mapsto \alpha_j(t,\log(\vert{y_0(t,x)}\vert)),$$
for $j \in \{1,...,m+s\}$ we are done by using composition of power series.

\textbf{Case 2:} $\vartheta \neq 0$. There is $\epsilon \in \textnormal{}]0,1[$ such that $0<\vert{w_0-\vartheta(z,w')}\vert < \epsilon \vert{w_0}\vert$ for $(z,w) \in D_A$. This gives with Remark \ref{2.12} $w_0 \sim_{D_A} \vartheta(z,w')$ and therefore
$$y_0 \sim_A \vartheta(\eta,y_1,...,y_r).$$
By Claim 2 $f_j|_A$ is a nicely log-analytic function in $x$ of order at most $r-1$ for every $j \in \{1,...,m\}$. With the inductive hypothesis applied to $f_j|_A$ for $j \in \{1,...,m\}$ we are done.

\hfill$\blacksquare$

Since every log-analytic function is nicely log-analytic in $x$, Theorem A is established.

{\bf Theorem A}\\
{\it
Let $m \in \mathbb{N}$, $r \in \mathbb{N}_0$. Let $X \subset \mathbb{R}^n \times \mathbb{R}$ be definable. Let $f_1,....,f_m:X \to \mathbb{R}$ be log-analytic functions of order at most $r$. Then there is a definable cell decomposition $\mathcal{C}$ of $X_{\neq 0}$ such that $f_1|_C,...,f_m|_C$ are nicely $r$-log-analytically prepared in $x$ in a simultaneous way for every $C \in \mathcal{C}$.}
 
Outgoing from Theorem A one may ask if there is an another kind of preparation of log-analytic functions of order $>0$ with log-analytic data only. To investigate this question one notes that an interesting example is the class of constructible functions introduced by Cluckers and Miller in \cite{10.1215/00127094-2010-213}. A function $f:\mathbb{R}^m \to \mathbb{R}$ is called constructible if it is a finite sum of finite products of functions on $\mathbb{R}^m$ which are definable in $\mathbb{R}_{\textnormal{an}}$ and of functions which are the logarithm of a positive function on $\mathbb{R}^m$ which is definable in $\mathbb{R}_{\textnormal{an}}$. This is a proper larger class than the globally subanalytic functions, but a proper subclass of the class of log-analytic ones of order at most $1$ (since the function $\mathbb{R} \to \mathbb{R}, x \mapsto \sqrt{\log(\vert{x}\vert+1)},$ is log-analytic of order $1$ but not constructible). For constructible functions there is a pure preparation not in terms of units but suitable for questions on integration (compare with \cite{10.1215/00127094-2010-213}).

A important consequence of Theorem A is the following: Call a definable cell $C \subset \mathbb{R}^n \times \mathbb{R}$ \textbf{simple} if $C_t$ is of the form $]0,d_t[$ for every $t \in \pi(C)$ where $d_t \in \mathbb{R}_{>0} \cup \{\infty\}$ (see for example Definition 2.15 in \cite{KaiserOpris}). Proposition 2.19 in \cite{KaiserOpris} states that the center of every logarithmic scale on such a simple cell $C$ vanishes. This gives the following (compare with Theorem 2.30 in \cite{KaiserOpris}).

\begin{cor}\label{2.55}
Let $X \subset \mathbb{R}^n \times \mathbb{R}$ be definable. Let $f:X \to \mathbb{R}$ be log-analytic of order at most $r$. Then there is a definable cell decomposition $\mathcal{C}$ of $X_{\neq 0}$ such that for every simple $C \in \mathcal{C}$ we have that $f|_C$ is purely $r$-log-analytically prepared in $x$ with center $0$.
\end{cor}

{\bf Proof:} By Theorem A there is a definable cell decomposition $\mathcal{C}$ of $X_{\neq 0}$ such that for every $C \in \mathcal{C}$ we have that $f|_C$ is nicely $r$-log-analytically prepared in $x$. Fix a simple $C \in \mathcal{C}$ and let 
$$(r,\mathcal{Y},a,q,s,v,b,P)$$
be a nice LA-preparing tuple for $f$. Let $\Theta$ be the center of $\mathcal{Y}$. Note that $\Theta=0$. Let $E$ be a set of $C$-heirs such that $a$ and $b$ can be constructed from $E$. Let $h \in E$. There is $\hat{r} \in \mathbb{N}_0$, an $\hat{r}$-logarithmic scale $\hat{\mathcal{Y}}$ with center $(\hat{\Theta}_0,...,\hat{\Theta}_{\hat{r}})$ on $C$ and $l \in \{1,...,\hat{r}\}$ such that $h=\exp(\hat{\Theta}_l)$. Note that $\hat{\Theta}_l=0$. So we have $h=1$. So we obtain $E=\emptyset$ or $E=\{1\}$. With the proof of the claim in (2) of Example \ref{2.32} one sees that $a$ and $b$ are log-analytic. \hfill$\blacksquare$

Important consequences of Corollary \ref{2.55} are differentiability results for the class of log-analytic functions like strong quasianalyticity and a parametric version of Tamm's Theorem (see \cite{KaiserOpris}).

\section{Preparation Theorems for Definable Functions in $\mathbb{R}_{\textnormal{an,exp}}$}

This section is devoted to the proofs of Theorem B and C, the latter being the main result of the paper. We start with some preparations which we need for the rigorous proof of this both theorems.

\subsection{Preparations}
For Section 3.1 we fix $m \in \mathbb{N}$, a tuple of variables $v:=(v_1,...,v_m)$, and a definable set $X \subset \mathbb{R}^m$. Fix $k,l \in \mathbb{N}_0$ and definable functions $g_1,...,g_k:X \to \mathbb{R}$ and $h_1,...,h_l:X \to \mathbb{R}$. Set $\beta:=(g,h_l,\exp(h_1),...,\exp(h_l))$. (We could also write $\beta=(g,\exp(h_1),...,\exp(h_l))$, but this simplifies notation below.) Note that $\beta(X) \subset \mathbb{R}^{k+l} \times \mathbb{R}_{>0}$. So there exists a $0$-logarithmic scale on $\beta(X)$, i.e. $\beta(X)$ is $0$-admissible. 

Fix a log-analytic function $F:\mathbb{R}^{k+l+1} \to \mathbb{R}$. Let $\alpha:X \to \mathbb{R}, v \mapsto F(\beta(v))$. Let $y:=(y_1,...,y_{k+l})$ range over $\mathbb{R}^{k+l}$. Let $z$ be another single variable such that $(y,z)$ ranges over $\mathbb{R}^{k+l} \times \mathbb{R}$.  Let $\pi^*:\mathbb{R}^{k+l} \times \mathbb{R} \to \mathbb{R}^{k+l}, (y,z) \mapsto y$. 

\begin{prop}\label{3.1}
Let $\Theta:\mathbb{R}^{k+l} \to \mathbb{R}$ be log-analytic such that 
$$\exp(h_l) \sim_X \Theta(g,h_l,\exp(h_1),...,\exp(h_{l-1})).$$ 
There is a log-analytic function $G:\mathbb{R}^{k+l} \to \mathbb{R}$ such that 
$$\alpha = G(g,h_l,\exp(h_1),...,\exp(h_{l-1}))$$
on $X$. 
\end{prop}

{\bf Proof:}
Let $\kappa:=\Theta(g,h_l,\exp(h_1),...,\exp(h_{l-1}))$. Note that $\kappa>0$. There is $\delta>1$ such that
$$\frac{1}{\delta} < \frac{\exp(h_l)}{\kappa} < \delta$$
on $X$. By taking logarithm we get with $\lambda:=\log(\delta)$
$$-\lambda < h_l-\log(\kappa) < \lambda$$
on $X$. Set 
$$
\exp^*:\mathbb{R} \to \mathbb{R},\textnormal{ } x \mapsto 
\left\{\begin{array}{ll} \exp(x), & x \in [-\lambda,\lambda], \\
0,& \textnormal{ otherwise.} \end{array}\right. 
$$
Then $\exp^*$ is globally subanalytic and we have 
$$\alpha= F(g,h_l,\exp(h_1),...,\exp(h_{l-1}),\kappa \cdot \exp^*(h_l-\log(\kappa)))$$
on $X$. Consider 
$$G:\mathbb{R}^{k+l} \to \mathbb{R}, y \mapsto \left\{\begin{array}{lll} F(y,\Theta(y) \cdot \exp^*(y_{k+1}-\log(\Theta(y)))), & y \in \pi^*(\beta(X)), \\
0,& \textnormal{ otherwise.} \end{array}\right.$$
Note that $G$ is well-defined, log-analytic and we obtain $$\alpha=G(g,h_l,\exp(h_1),...,\exp(h_{l-1}))$$
on $X$ since for $x \in X$ we have 
$$(g(x),h_l(x),\exp(h_1(x)),...,\exp(h_{l-1}(x))) \in \pi^*(\beta(X)).$$
\hfill $\blacksquare$

\begin{prop}\label{3.2}
Assume that $F$ is positive and purely $0$-log-analytically prepared in $z$ with center $0$ on $\beta(X)$. Then there is a purely log-analytic function $H:\mathbb{R}^{k+l} \times \mathbb{R} \to \mathbb{R}$ in $z$ of order $0$ such that 
	$$\log(F(\beta))=H(\beta)$$
	on $X$.
\end{prop}

{\bf Proof:}
Let
$$(0,\mathcal{Y},a,q,s,v,b,P)$$
be a pure LA-preparing tuple for $f$ where $b:=(b_1,...,b_s)$ and $\mathcal{Y}:=y$. Note that $a>0$. Consider
$$\eta:=(\eta_0,...,\eta_s):\pi^*(\beta(X)) \to \mathbb{R}_{>0} \times \mathbb{R}^s,y \mapsto (a(y),b_1(y),...,b_s(y)).$$
Note that $\eta$ is log-analytic. Let $w:=(w_0,...,w_{s+2})$ range over $\mathbb{R}^{s+3}$. For $w \in \mathbb{R}_{>0} \times \mathbb{R}^s \times \mathbb{R} \times \mathbb{R}_{\neq 0}$ with $-1 \leq w_i\vert{w_{s+2}}\vert^{p_{0i}} \leq 1$ for every $i \in \{1,...,s\}$ let 
$$\phi(w):=\log(w_0)+qw_{s+1}+\log(v(w_1 \vert{w_{s+2}}\vert^{p_{01}},...,w_s \vert{w_{s+2}}\vert^{p_{0s}})).$$
Consider 
$$G:\mathbb{R}_{>0} \times \mathbb{R}^s \times \mathbb{R} \times \mathbb{R}_{\neq 0} \to \mathbb{R}, w \mapsto$$ $$\left\{\begin{array}{ll} \phi(w), & -1 \leq w_i \vert{w_{s+2}}\vert^{p_{0i}} \leq 1 \textnormal{ for every } i \in \{1,...,s\},  \\
0, & \textnormal{otherwise.} \end{array}\right.$$
Then $G$ is purely log-analytic in $w_{s+2}$ of order $0$ since $\log(v)$ is globally subanalytic. Note that
$$\log(F(\beta)) = G(\eta(g,h_l,\exp(h_1),...,\exp(h_{l-1})),h_l,\exp(h_l))$$
on $X$. Then
$$
H:\mathbb{R}^{k+l} \times \mathbb{R} \to \mathbb{R}, (y,z) \mapsto \left\{\begin{array}{ll} G(\eta(y_1,...,y_{k+l}),y_{k+1},z) , & (y,z) \in \beta(X),  \\
0, & \textnormal{otherwise}, \end{array}\right.
$$
does the job, because $H$ is purely log-analytic in $z$ of order $0$. 
\hfill $\blacksquare$

\begin{cor}\label{3.3}
Let $c,d \in \mathbb{N}$. Suppose there are functions $\mu_1,...,\mu_c:\mathbb{R}^{k+l} \times \mathbb{R} \to \mathbb{R}$ and $\nu_1,...,\nu_d:\mathbb{R}^{k+l} \times \mathbb{R} \to \mathbb{R}_{>0}$ which are purely $0$-log-analytically prepared in $z$ with center $\Theta$ on $\beta(X)$, and a globally subanalytic function $G:\mathbb{R}^{c+d} \to \mathbb{R}$ such that    
	$$\alpha=G(\mu_1(\beta),...,\mu_c(\beta),\log(\nu_1(\beta)),...,\log(\nu_d(\beta))).$$
	If $\Theta=0$ there is a purely log-analytic function $H:\mathbb{R}^{k+l} \times \mathbb{R} \to \mathbb{R}$ in $z$ of order $0$ such that $\alpha=H(\beta)$ on $X$. If $\Theta \neq 0$ then we have
	$$\exp(h_l) \sim_X \Theta(g,h_l,\exp(h_1),...,\exp(h_{l-1})).$$
\end{cor}

{\bf Proof:}
Note that $y_0:\beta(X) \to \mathbb{R},(y,z) \mapsto z-\Theta(y)$, is a $0$-logarithmic scale with log-analytic center $\Theta$.

{\bf Case 1:} $\Theta = 0$. Then by Proposition \ref{3.2} there are purely log-analytic functions $H_1,...,H_m:\mathbb{R}^{k+l} \times \mathbb{R} \to \mathbb{R}$ in $z$ of order $0$ such that $\log(\nu_j(\beta))=H_j(\beta)$ on $X$ for every $j \in \{1,...,m\}$. By (2) in Remark \ref{2.27} we see that $\mu_1,...,\mu_c$ are purely log-analytic in $z$ of order $0$. Because $G$ is globally subanalytic we are done.

{\bf Case 2:} $\Theta \neq 0$. By Remark \ref{2.12} it is $z \sim_{\beta(X)} \Theta$. We obtain the result.
\hfill$\blacksquare$

\begin{cor}\label{3.4}
	Let $r \in \mathbb{N}_0$. Suppose $F$ is purely log-analytic in $z$ of order $r$. Then there is a decomposition $\mathcal{C}$ of $X$ into finitely many definable cells such that for $C \in \mathcal{C}$ the following holds. There is a purely log-analytic function $H:\mathbb{R}^{k+l} \times \mathbb{R} \to \mathbb{R}$ in $z$ of order $0$ such that $\alpha|_C=H(\beta|_C)$ or there is a log-analytic function $\Theta:\mathbb{R}^{k+l} \to \mathbb{R}$ such that
	$$\exp(h_l) \sim_C \Theta(g,h_l,\exp(h_1),...,\exp(h_{l-1})).$$
\end{cor}

{\bf Proof:}
We do an induction on $r$.\\
$r=0$: Then $F$ is purely log-analytic in $z$ of order $0$. The assertion follows.\\
$r-1 \to r$: With Definition \ref{2.25} it is enough to consider the following situation. Let $c,d \in \mathbb{N}_0$, a globally subanalytic function $G:\mathbb{R}^{c+d} \to \mathbb{R}$, purely log-analytic functions $\rho_1,...,\rho_c:\mathbb{R}^{k+l} \times \mathbb{R} \to \mathbb{R}$ and $\sigma_1,...,\sigma_d:\mathbb{R}^{k+l} \times \mathbb{R} \to \mathbb{R}_{>0}$ in $z$ of order at most $r-1$ be such that  
$$\alpha=G(\rho_1(\beta),...,\rho_c(\beta),\log(\sigma_1(\beta)),...,\log(\sigma_d(\beta))).$$
Applying the inductive hypothesis to $\rho_1(\beta),...,\rho_c(\beta),\sigma_1(\beta),...,\sigma_d(\beta)$ there is a decomposition $\mathcal{A}$ of $X$ into finitely many definable cells such that the following holds for $A \in \mathcal{A}$. There are functions $\mu_1,...,\mu_c:\mathbb{R}^{k+l} \times \mathbb{R} \to \mathbb{R}$ and $\nu_1,...,\nu_d:\mathbb{R}^{k+l} \times \mathbb{R} \to \mathbb{R}_{>0}$ which are purely log-analytic in $z$ of order $0$ such that
$$\alpha=G(\mu_1(\beta),...,\mu_c(\beta),\log(\nu_1(\beta)),...,\log(\nu_d(\beta)))$$
on $A$ ($G$ does not change) or there is a log-analytic function $\tilde{\Theta}:\mathbb{R}^{k+l} \to \mathbb{R}$ such that
$$\exp(h_l) \sim_A \tilde{\Theta}(g,h_l,\exp(h_1),...,\exp(h_{l-1})).$$
Fix such an $A$. If the latter holds we are done. So assume the former. Fix the corresponding $\mu_1,...,\mu_c,\nu_1,...,\nu_d$. By Remark \ref{2.53} there is a definable cell decomposition $\mathcal{D}$ of $\mathbb{R}^{k+l} \times \mathbb{R}_{\neq 0}$ such that $\mu_1,...,\mu_c,\nu_1,...,\nu_d$ are purely $0$-log-analytically prepared in $z$ in a simultaneous way on every $D \in \mathcal{D}$. There is a definable cell decomposition $\mathcal{C}$ of $A$ such that for every $C \in \mathcal{C}$ there is $D_C \in \mathcal{D}$ with $\beta(C) \subset D_C$. Fix $C \in \mathcal{C}$ and the center $\Theta$ of the pure $0$-preparation of $\mu_1,...,\mu_c,\nu_1,...,\nu_d$ on $D_C$. If $\Theta=0$ we find with Corollary \ref{3.3} a purely log-analytic function $H:\mathbb{R}^{k+l} \times \mathbb{R} \to \mathbb{R}$ in $z$ of order $0$ such that $\alpha|_C = H(\beta|_C)$. If $\Theta \neq 0$ then with Remark \ref{2.12} $z \sim_{D_C} \Theta(y)$ and therefore
\[
\exp(h_l) \sim_C \Theta(g,h_l,\exp(h_1),...,\exp(h_{l-1})).
\]
\hfill$\blacksquare$
\subsection{Proof of Theorem B and Theorem C}
For a definable function $f:X \to \mathbb{R}$ (where $X \subset \mathbb{R}^m$ is definable) there exist $e \in \mathbb{N}_0$ and a set $E$ of positive definable functions on $X$ such that $f$ has exponential number at most $e$ with respect to $E$ (i.e. is a composition of log-analytic functions and exponentials from $E$, see Definition \ref{1.5} and Remark \ref{1.7}). \\
The first goal for this section is to prepare $f$ cellwise as a product of a log-analytic function, an exponential of a function $h$ which is the restriction of a finite $\mathbb{Q}$-linear combination of functions from $\log(E)$ with exponential number at most $e-1$ with respect to $E$, and a unit of a special form (see Theorem B). Additionally $h$ is itself prepared. To be more precise we introduce cellwise a set $P$ of positive definable functions (or exponentials) which occur in the preparation of $f$ and show that every function from $\log(P)$ is also prepared (i.e. a product of a log-analytic function, an exponential of a function $g$ with $\exp(g) \in P$ and a unit of a special form). So analytical properties of functions from $\log(E)$ closed under taking finite $\mathbb{Q}$-linear combinations (like the property of being locally bounded) can be cellwise transferred to every function in $\log(P)$ (see \cite{10.1215/00192082-10234104}).\\ 
The second goal for this section is to prove Theorem C by combining Theorem A and Theorem B: From Theorem B we also get cellwise a finite set $L$ of log-analytic functions which occur in the preparation of $f$ and with Theorem A we prepare all functions from $L$ simultaneously.

\begin{defn}\label{3.5}
Let $X \subset \mathbb{R}^m$ be definable and $f:X \to \mathbb{R}$ be a function. Let $E$ be a set of positive definable functions on $X$.

\begin{itemize}
	\item [(1)]
	Let $L$ be a set of log-analytic functions on $X$. By induction on $e \in \mathbb{N}_0 \cup \{-1\}$ we define that $f:X \to \mathbb{R}$ is \textbf{$e$-prepared with respect to $L$ and $E$}.
	
	\textbf{Base Case}: The function $f$ is $(-1)$-prepared with respect to $L$ and $E$ if $f$ is the zero function.
	
	\textbf{Inductive step}: The function $f$ is $e$-prepared with respect to $L$ and $E$ if 
	$$f=a \cdot \exp(c) \cdot u$$
	where $a\in L$ vanishes or does not have a zero, $c:X \to \mathbb{R}$ is $(e-1)$-prepared with respect to $L$ and $E$, $\exp(c) \in E$ and $u:X \to \mathbb{R}$ is a function of the form
	$$u=v(b_1 \cdot \exp(d_1),...,b_s \cdot \exp(d_s))$$
	where $s \in \mathbb{N}_0$, $b_j \in L$ does not have any zero, $d_j:X \to \mathbb{R}$ is $(e-1)$-prepared with respect to $L$ and $E$ and $\exp(d_j) \in E$ for every $j \in \{1,...,s\}$.
	
	Additionally $v$ is a power series which converges absolutely on an open neighborhood of $[-1,1]^s$, it is $b_j(x)\exp(d_j(x)) \in [-1,1]$ for every $x \in X$ and every $j \in \{1,...,s\}$, and $v([-1,1]^s) \subset \mathbb{R}_{>0}$.
	
	\item [(2)] We say that $f$ is \textbf{$e$-prepared with respect to $E$} if there is a set $L$ of log-analytic functions on $X$ such that $f$ is $e$-prepared with respect to $L$ and $E$.  
\end{itemize}
\end{defn}

\begin{rem}\label{3.6}
Let $X \subset \mathbb{R}^m$ be definable. The following holds.
\begin{itemize}
	\item [(1)]
	Let $e \geq 0$. If $f:X \to \mathbb{R}$ is $e$-prepared with respect to a set $E$ of positive definable functions then $1 \in E$.
	\item [(2)]
	Let $f:X \to \mathbb{R}$ be log-analytic. Then $f$ is $0$-prepared with respect to $L:=\{f\}$ and $E:=\{1\}$.
	\item [(3)]
	If $f$ is $0$-prepared with respect to a set $E$ of positive definable functions then $f$ is log-analytic.
\end{itemize}
\end{rem}

\begin{rem}\label{3.7}
Let $X \subset \mathbb{R}^m$ be definable. Let $E$ be a set of positive definable functions and $L$ be a set of log-analytic functions on $X$. Let $e \in \mathbb{N}_0$. Let $f:X \to \mathbb{R}$ be $e$-prepared with respect to $L$ and $E$. There are $k \in \mathbb{N}$, log-analytic functions $h_1,...,h_k \in L$, $(e-1)$-prepared $g_1,...,g_k$ with respect to $L$ and $E$ with $\exp(g_1),...,\exp(g_k) \in E$ and a globally subanalytic function $G:\mathbb{R}^{2k} \to \mathbb{R}$ such that
$$f=G(h_1,...,h_k,\exp(g_1),...,\exp(g_k)).$$
\end{rem}

{\bf Proof:}
The proof is similar as the proof of (1) in Proposition \ref{2.43}. \hfill$\blacksquare$

\iffalse
There are $s \in \mathbb{N}$, $a,b_1,...,b_s \in L$ and $(e-1)$-prepared $c,d_1,....,d_s$ with respect to $L$ and $E$ such that $\exp(c),\exp(d_1),...,\exp(d_s) \in E$ and
$$f=a \exp(c)v(b_1\exp(d_1),...,b_s\exp(d_s))$$
where $v$ is a power series which converges absolutely on an open neighborhood of $[-1,1]^s$ and $b_j(x)\exp(d_j(x)) \in \textnormal{}[-1,1]$ for every $j \in \{1,...,s\}$ and every $x \in X$. Note that the function 
$$H:\mathbb{R}^s \to \mathbb{R}, (u_1,...,u_s) \mapsto \left\{\begin{array}{ll} v(u_1,...,u_s), & (u_1,...,u_s) \in [-1,1]^s, \\
0, & \textnormal{otherwise}, \end{array}\right. $$
is globally subanalytic. Choose $k:=s+1$, $h_1:=a$, $h_j:=b_{j-1}$ for $j \in \{2,...,k\}$, $g_1:=c$ and $g_j:=d_{j-1}$ for $j \in \{2,...,k\}$. Let $w:=(w_1,...,w_k)$ and $z:=(z_1,...,z_k)$ range over $\mathbb{R}^k$. Consider 
$$G:\mathbb{R}^{2k} \to \mathbb{R}, (w,z) \mapsto w_1z_1H(w_2z_2,...,w_kz_k).$$
Then $G$ is globally subanalytic and we have 
$$f=G(h_1,...,h_k,\exp(g_1),...,\exp(g_k)).$$

\hfill$\blacksquare$
\fi

\begin{rem}\label{3.8}
Let $X \subset \mathbb{R}^m$ be definable. Let $E$ be a set of positive definable functions on $X$. Let $e \in \mathbb{N}_0$. Let $f:X \to \mathbb{R}$ be $e$-prepared with respect to $E$. Then $f$ has exponential number at most $e$ with respect to $E$ and is therefore definable.
\end{rem}

{\bf Proof:}
This is easily seen with Remark \ref{3.7} and an easy induction on $e$. \hfill$\blacksquare$
\iffalse
We do an induction on $e$. If $e=0$ then $f$ is log-analytic by (3) in Remark \ref{3.6}. So assume $e>0$. Let $L$ be a set of log-analytic functions such that $f$ is $e$-prepared with respect to $L$ and $E$. By Remark \ref{3.6} there are $k \in \mathbb{N}$, log-analytic functions $h_1,...,h_k \in L$, $(e-1)$-prepared $g_1,...,g_k$ with respect to $L$ and $E$ with $\exp(g_1),...,\exp(g_k) \in E$ and a globally subanalytic function $G:\mathbb{R}^{2k} \to \mathbb{R}$ such that
$$f=G(h_1,...,h_k,\exp(g_1),...,\exp(g_k)).$$
With the inductive hypothesis and Definition \ref{1.5}(a) we are done.
\fi

Now we are ready to formulate and prove Theorem B. 

{\bf Theorem B}\\
{\it Let $X \subset \mathbb{R}^m$ be definable and $e \in \mathbb{N}_0$. Let $f:X \to \mathbb{R}$ be a function. Let $E$ be a set of positive definable functions on $X$ such that $f$ has exponential number at most $e$ with respect to $E$. Then there is a decomposition $\mathcal{C}$ of $X$ into finitely many definable cells such that for every $C \in \mathcal{C}$ there is a finite set $P$ of positive definable functions on $C$ and a finite set $L$ of log-analytic functions on $C$ such that the following holds.
\begin{itemize}
		\item [(1)]
		$f|_C$ is $e$-prepared with respect to $L$ and $P$ and for every $g \in \log(P)$ there is $\alpha \in \{-1,...,e-1\}$ such that $g$ is $\alpha$-prepared with respect to $L$ and $P$.
		\item [(2)]
		$P$ satisfies the following condition $(*_e)$ with respect to $E|_C$: If $g \in \log(P)$ is $l$-prepared with respect to $L$ and $P$ for $l \in \{0,...,e-1\}$ then $g$ is a finite $\mathbb{Q}$-linear combination of functions from $\log(E)$ restricted to $C$ which have exponential number at most $l$ with respect to $E$.
\end{itemize}}

{\bf Proof:}
We do an induction on $e$. \\
$e=0:$ Then $f$ is log-analytic and we are done by choosing $P=\{1\}$ and $L=\{f\}$.\\
$e-1 \to e:$ There are $k,l \in \mathbb{N}$, functions $g_1,...,g_k,h_1,...,h_l:\mathbb{R}^m \to \mathbb{R}$ with exponential number at most $e-1$ with respect to $E$, and a log-analytic function $F:\mathbb{R}^{k+l} \to \mathbb{R}$ such that
$$f = F(g_1,...,g_k,\exp(h_1),...,\exp(h_l))$$
and $\exp(h_1),...,\exp(h_l) \in E$. 

By an auxiliary induction on $l$ and involving the inductive hypothesis we may assume that the theorem is proven for functions of the form 
$$\kappa = H(g_1,...,g_k,h_l,\exp(h_1),...,\exp(h_{l-1}))$$
on $X$ where $H:\mathbb{R}^{k+l} \to \mathbb{R}$ is log-analytic. $(**)$

(If $l=1$ then $(**)$ holds by the inductive hypothesis: $g = H(g_1,...,g_k,h_l)$ has exponential number at most $e-1$ with respect to $E$ by Remark \ref{1.8}.)

(This includes also functions of the form 
$$\kappa = \hat{H}(g_1,...,g_k,\exp(h_1),...,\exp(h_{l-1}))$$
on $X$ where $\hat{H}:\mathbb{R}^{k+l-1} \to \mathbb{R}$ is log-analytic.)

Let $g:=(g_1,...,g_k)$. Let $y:=(y_1,...,y_{k+l})$ range over $\mathbb{R}^{k+l}$. Let $z$ be another single variable such that $(y,z)$ ranges over $\mathbb{R}^{k+l} \times \mathbb{R}$. Let $y':=(y_1,...,y_{k+l-1})$. Let $\pi^+:\mathbb{R}^{k+l} \to \mathbb{R}^{k+l-1},y \mapsto y'$. Let $r \in \mathbb{N}_0$ be such that $F$ is purely log-analytic in $y_{k+l}$ of order at most $r$.

\textbf{Case 1}: $r=0$. By Remark \ref{2.53} we find a definable cell decomposition $\mathcal{D}$ of $\mathbb{R}^{k+l-1} \times \mathbb{R}_{\neq 0}$ such that $F|_D$ is purely $0$-log-analytically prepared in $y_{k+l}$ for every $D \in \mathcal{D}$. There is a decomposition $\mathcal{A}$ of $X$ into finitely many definable cells such that for every $A \in \mathcal{A}$ there is $D_A \in \mathcal{D}$ such that for every $x \in A$ we have $(g(x),\exp(h_1(x)),...,\exp(h_l(x))) \in D_A$. %, because $\exp(h_l)>0$ on $C$. 
Fix $A \in \mathcal{A}$ and a purely preparing tuple $(0,\mathcal{Y},a,q,s,v,b,P)$ for $F|_{D_A}$ where $b:=(b_1,...,b_s)$ and $P:=(p_1,...,p_s)^t$. Let $\Theta$ be the center of $\mathcal{Y}$. Then one of the following properties holds.

\begin{itemize}
	\item [(1)] There is $\epsilon \in \textnormal{}]0,1[$ such that
	$$0< \vert{y_{k+l} - \Theta(y')}\vert <\epsilon \vert{y_{k+l}}\vert$$
	for all $y \in D_A$.
	\item [(2)] %$w_k \neq 0$ for all $(z,w) \in E_D$. 
	$\Theta=0$ on $\pi^+(D_A)$.
\end{itemize}

Assume $(1)$. Then by Remark \ref{2.12} we have $y_{k+l} \sim_{D_A} \Theta$. This gives 
$$\exp(h_l) \sim_A \Theta(g,\exp(h_1),...,\exp(h_{l-1})).$$
By Proposition \ref{3.1} there is a log-analytic function $G:\mathbb{R}^{k+l} \to \mathbb{R}$ such that
$$f=G(g,h_l,\exp(h_1),...,\exp(h_{l-1}))$$
on $A$. With $(**)$ applied to $f$ we are done.

Assume $(2)$. Let $\beta:=(g,\exp(h_1),...,\exp(h_{l-1}))$. It holds
$$f|_A=a(\beta)\exp(qh_l)v(b_1(\beta)\exp(p_1h_l),...,b_s(\beta)\exp(p_sh_l)).$$
Let $b_0:=a$. Note that we can apply $(**)$ to $b_j(\beta)$ for $j \in \{0,...,s\}$ and so there is a decomposition $\mathcal{B}$ of $A$ into finitely many definable cells such that for every $B \in \mathcal{B}$ the following holds: There is a finite set $P'$ of positive definable functions on $B$ and a finite set $L'$ of log-analytic functions on $B$ such that $b_0(\beta)|_B,...,b_s(\beta)|_B$ are $e$-prepared with respect to $L'$ and $P'$ and that $P'$ satisfies property $(*_e)$ with respect to $E|_B$. Fix such a $B$.
Then we have for $j \in \{0,...,s\}$
$$b_j(\beta)=\hat{a}_j\exp(\hat{c}_j)\hat{v}_j(\hat{b}_{1j}\exp(\hat{d}_{1j}),...,\hat{b}_{s_jj}\exp(\hat{d}_{s_jj}))$$
where $s_j \in \mathbb{N}$, $\hat{a}_j,\hat{b}_{1j},...,\hat{b}_{s_jj}$ are log-analytic and $\hat{c}_j,\hat{d}_{1j},...,\hat{d}_{s_jj}$ are finite $\mathbb{Q}$-linear combinations of functions from $\log(E)$ restricted to $B$ which have exponential number at most $e-1$ with respect to $E$ and $\hat{v}_j$ is a power series which converges absolutely on an open neighborhood of $[-1,1]^{s_j}$ and $\hat{v}_j([-1,1]^{s_j}) \subset \mathbb{R}_{>0}$. Since $h_l \in \log(E)$ has exponential number at most $e-1$ with respect to $E$ we obtain that $\hat{c}_0+qh_l$ and $\hat{c}_j+p_jh_l$ for $j \in \{1,...,s\}$ are finite $\mathbb{Q}$-linear combinations of functions from $\log(E)$ restricted to $B$ which have exponential number at most $e-1$ with respect to $E$. (Note that $h_l$ is not necessarily $(e-1)$-prepared with respect to $P'$ on $B$.)

Consequently we obtain by composition of power series that there are $\tilde{s} \in \mathbb{N}$, log-analytic functions $\tilde{a},\tilde{b}_1,...,\tilde{b}_{\tilde{s}}:B \to \mathbb{R}$, and functions $c,d_1,...,d_{\tilde{s}}:B \to \mathbb{R}$ which are finite $\mathbb{Q}$-linear combinations of elements from $\log(E)$ restricted to $B$ which have exponential number at most $e-1$ with respect to $E|_B$ such that the following holds: we have
$$f|_B=\tilde{a}\exp(c) \tilde{v}(\tilde{b}_1\exp(d_1),...,\tilde{b}_{\tilde{s}}\exp(d_{\tilde{s}}))$$
where $\tilde{b}_j(x)\exp(d_j(x)) \in [-1,1]$ for every $x \in B$ and $j \in \{1,...,\tilde{s}\}$, and $\tilde{v}$ is a power series which converges absolutely on an open neighborhood of $[-1,1]^{\tilde{s}}$ and $\tilde{v}([-1,1]^{\tilde{s}}) \subset \mathbb{R}_{>0}$. Note that $c$ and $d_1,...,d_{\tilde{s}}$ have exponential number at most $e-1$ with respect to $E$. By the inductive hypothesis there is a further decomposition $\mathcal{C}_B$ of $B$ into finitely many definable cells such that for every $C \in \mathcal{C}_B$ there is a finite set $\tilde{P}$ of positive definable functions on $C$ which satisfies property $(*_{e-1})$ with respect to $E|_C$ and a finite set $\tilde{L}$ of log-analytic functions on $C$ such that the functions $c|_C$ and $d_1|_C,...,d_{\tilde{s}}|_C$ are $(e-1)$-prepared with respect to $\tilde{L}$ and $\tilde{P}$. So choose 
$$L:=\tilde{L} \cup \{\tilde{a}|_C,\tilde{b}_1|_C,...,\tilde{b}_{\tilde{s}}|_C\}$$
and
$$P:=\tilde{P} \cup \{\exp(c)|_C,\exp(d_1)|_C,...,\exp(d_{\tilde{s}})|_C\}.$$
Then $P$ satisfies property $(*_e)$ with respect to $E|_C$. Note that $f|_C$ is $e$-prepared with respect to $L$ and $P$ and we are done.

\textbf{Case 2}: $r>0$. By Corollary \ref{3.4} there is a decomposition $\mathcal{A}$ of $X$ into finitely many definable cells such that for every $A \in \mathcal{A}$ one of the following properties holds.
\begin{itemize}
	\item [(1)]
	There is a purely log-analytic function $H:\mathbb{R}^{k+l} \times \mathbb{R} \to \mathbb{R}$ in $z$ of order $0$ such that
	$$f|_A=H(g,h_l,\exp(h_1),...,\exp(h_l)).$$
	\item [(2)]
	There is a log-analytic function $\tilde{\Theta}:\mathbb{R}^{k+l} \to \mathbb{R}$ such that 
	$$\exp(h_l) \sim_A \tilde{\Theta}(g,h_l,\exp(h_1),...,\exp(h_{l-1})).$$ 
\end{itemize}
Let $A \in \mathcal{A}$. If $(1)$ holds for $A$ then we are done with Case 1. If $(2)$ holds for $A$ then by Proposition \ref{3.1} there is a log-analytic function $H:\mathbb{R}^{k+l} \to \mathbb{R}$ such that 
$$f=H(g,h_l,\exp(h_1),...,\exp(h_{l-1}))$$
on $A$. We are done with $(**)$ applied to $f$.
\hfill$\blacksquare$

In the following we aim for the second goal of Section 3.2, the proof of Theorem C. In contrast to the treatment of Lion-Rolin in \cite{Lion1997} we describe precisely how the unit in the preparation in Theorem C looks like, bound the number of iterations of the exponentials which occur in the preparation by the exponential number of $f$ and determine which exponential functions occur: From Theorem B we cellwise take over the set $P$ of exponentials (each function from $\log(P)$ is a restriction of a finite $\mathbb{Q}$-linear combination of functions from $\log(E)$) and Theorem A, the log-analytic preparation theorem, gives $C$-nice data (i.e. $C$-heirs). \\
Furthermore the preparation of every $g \in \log(P) \cup \{f\}$ is described by a pair $(l,r) \in (\mathbb{N}_0 \cup \{-1\}) \times \mathbb{N}_0$: From Theorem $B$ we obtain finite sets $L$,$P$ and an $l \in \mathbb{N}_0 \cup \{-1\}$ such that $g$ is $l$-prepared with respect to $L$ and $P$ and from Theorem A an $r \in \mathbb{N}_0$ such that the functions from $L$ are nicely $r$-log-analytically prepared in a simultaneous way (the number $r$ does not depend on $g$ and is chosen so that every function from $L$ is log-analytic of order at most $r$). 

After the complete proof of Theorem C we will give some important consequences.

For the rest of this section let $n \in \mathbb{N}_0$, $t:=(t_1,...,t_n)$ range over $\mathbb{R}^n$, $x$ over $\mathbb{R}$ and let $\pi:\mathbb{R}^{n+1} \to \mathbb{R}^n, (t,x) \mapsto t$.

\begin{defn}\label{3.9}
Let $r \in \mathbb{N}_0$ and $e \in \mathbb{N}_0 \cup \{-1\}$. Let $C \subset \mathbb{R}^n \times \mathbb{R}$ be definable and $f:C \to \mathbb{R}$ be a function. Let $E$ be a set of positive definable functions on $C$. By induction on $e$ we define that $f$ is \textbf{$(e,r)$-prepared} in $x$ with center $\Theta$ with respect to $E$. To this preparation we assign a \textbf{preparing tuple} for $f$.

$e=-1$: We say that $f$ is $(-1,r)$-prepared in $x$ with center $\Theta$ with respect to $E$ if $f$ is the zero function. A preparing tuple for $f$ is $(0)$.

$e-1 \to e$: We say that $f$ is $(e,r)$-prepared in $x$ with center $\Theta$ with respect to $E$ if 
$$f(t,x)=a(t)\vert{\mathcal{Y}(t,x)}\vert^{\otimes q}e^{c(t,x)} u(t,x)$$
for every $(t,x) \in C$ where $a:\pi(C) \to \mathbb{R}$ is $C$-nice which vanishes identically or has no zero, $\mathcal{Y}:=(y_0,...,y_r)$ is a nice $r$-logarithmic scale with center $\Theta$, $q=(q_0,...,q_r) \in \mathbb{Q}^{r+1}$, $\exp(c) \in E$ where $c$ is $(e-1,r)$-prepared in $x$ with center $\Theta$ with respect to $E$ and $u=v \circ \phi$ where $v$ is a power series which converges on an open neighborhood of $[-1,1]^s$ with $v([-1,1]^s) \subset \mathbb{R}_{>0}$ and $\phi:=(\phi_1,...,\phi_s):C \to [-1,1]^s$ is a function of the form 
$$\phi_j:=b_j(t)\vert{y_0}\vert^{p_{j0}} \cdot ... \cdot \vert{y_r}\vert^{p_{jr}}\exp(c_j(t,x))$$
for $j \in \{1,...,s\}$ where $b_j:\pi(C) \to \mathbb{R}$ is $C$-nice, $p_{j0},...,p_{jr} \in \mathbb{Q}$, $\exp(c_j) \in E$ and $c_j$ is $(e-1,r)$-prepared in $x$ with center $\Theta$ with respect to $E$. A preparing tuple for $f$ is then
$$\mathcal{J}:=(r,\mathcal{Y},a,\exp(c),q,s,v,b,\exp(d),P)$$
where $b:=(b_1,...,b_s)$, $\exp(d):=(\exp(d_1),...,\exp(d_s))$ and 
$$P:=\left(\begin{array}{cccc}
p_{10}&\cdot&\cdot&p_{1r}\\
\cdot&& &\cdot\\
\cdot&& &\cdot\\
p_{s0}&\cdot&\cdot&p_{sr}\\
\end{array}\right)\in M\big(s\times (r+1),\mathbb{Q}).$$
\end{defn}

\begin{rem}\label{3.10}
Let $e \in \mathbb{N}_0 \cup \{-1\}$ and $r \in \mathbb{N}_0$. Let $E$ be a set of positive definable functions on $C$. Let $f:C \to \mathbb{R}$ be $(e,r)$-prepared in $x$ with respect to $E$. If $e=0$ then $f$ is nicely log-analytically prepared in $x$.  
\end{rem}
\begin{rem}\label{3.11}
Let $e \in \mathbb{N}_0 \cup \{-1\}$ and $r \in \mathbb{N}_0$. Let $E$ be a set of positive definable functions on $C$. Let $f:C \to \mathbb{R}$ be $(e,r)$-prepared in $x$ with respect to $E$. Then there is a set $\mathcal{E}$ of $C$-heirs such that $f$ can be constructed from $E \cup \mathcal{E}$.
\end{rem}

{\bf Proof:}
We do an induction on $e$. For $e=-1$ there is nothing to show.\\
$e-1 \to e$: Let 
$$(r,\mathcal{Y},a,\exp(c),q,s,v,b,\exp(d),P)$$
be a preparing tuple for $f$ where $b:=(b_1,...,b_s)$, 
$$\exp(d):=(\exp(d_1),...,\exp(d_s))$$
and $\exp(c),\exp(d_1),...,\exp(d_s) \in E$. By the inductive hypothesis there is a set $\mathcal{E}'$ of $C$-heirs such that 
$c$ and $d_1,...,d_s$ can be constructed from $E \cup \mathcal{E}'$. By Remark \ref{1.8} we see that 
$\exp(c),\exp(d_1),...,\exp(d_s)$ can be constructed from $E \cup \mathcal{E}'$. Because $a$, $b_1,...,b_s$ and $\Theta_0,...,\Theta_r$ are $C$-nice there is a set $\bar{\mathcal{E}}$ of $C$-heirs such that $a$, $b_1,...,b_s$ and $\Theta_0,...,\Theta_r$ can be constructed from $\bar{\mathcal{E}}$. Set $\mathcal{E}:=\bar{\mathcal{E}} \cup \mathcal{E}'$. Then
$$\eta:=(a,\exp(c),b_1,...,b_s,\exp(d_1),...,\exp(d_s))$$
can be constructed from $E \cup \mathcal{E}$. Note also that $y_0,...,y_r$ can be constructed from $\mathcal{E}$ (compare with the proof of Remark \ref{2.48}). Let $k:=2+2s+r+1$. With a similar argument as in the proof of (1) in Proposition \ref{2.43} or Remark \ref{3.7} we see that there is a globally subanalytic $F:\mathbb{R}^k \to \mathbb{R}$ such that
$$f(t,x)=F(\eta(t),y_0(t,x),...,y_r(t,x)).$$
for every $(t,x) \in C$. With Remark \ref{1.8} we are done. \hfill$\blacksquare$

Now we are ready to formulate and prove Theorem C. 

{\bf Theorem C}\\
{\it Let $e \in \mathbb{N}_0$. Let $X \subset \mathbb{R}^n \times \mathbb{R}$ be definable and let $E$ be a set of positive definable functions on $X$. Let $f:X \to \mathbb{R}$ be a function with exponential number at most $e$ with respect to $E$. Then there is $r \in \mathbb{N}_0$ and a definable cell decomposition $\mathcal{C}$ of $X_{\neq 0}$ such that for every $C \in \mathcal{C}$ there is a finite set $P$ of positive definable functions on $C$ and $\Theta:=(\Theta_0,...,\Theta_r)$ such that the function $f|_C$ is $(e,r)$-prepared in $x$ with center $\Theta$ with respect to $P$. Additionally the following holds.
	\begin{itemize}
		\item [(1)]
		For every $g \in \log(P)$ there is $m \in \{-1,...,e-1\}$ such that $g$ is $(m,r)$-prepared in $x$ with center $\Theta$ with respect to $P$.
		\item [(2)]
		The following condition $(+_e)$ is satisfied: If $g \in \log(P)$ is $(l,r)$-prepared in $x$ with center $\Theta$ with respect to $P$ for $l \in \{-1,...,e-1\}$ then $g$ is a finite $\mathbb{Q}$-linear combination of functions from $\log(E)$ restricted to $C$ which have exponential number at most $l$ with respect to $E$.
\end{itemize}}

{\bf Proof:}
By Theorem B we may assume that $f$ is $e$-prepared with respect to a finite set $L:=\{g_1,...,g_m\}$ of log-analytic functions and a finite set $Q$ of positive definable functions with the following property: 

Every $g \in \log(Q)$ is $l$-prepared with respect to $L$ and $Q$ for an $l \in \{-1,...,e-1\}$. Additionally if $g \in \log(Q)$ is $l$-prepared for $l \in \{-1,...,e-1\}$ with respect to $L$ and $Q$ then $g$ is a finite $\mathbb{Q}$-linear combination of functions from $\log(E)$ which have exponential number at most $l$ with respect to $E$. $(*_e)$

Let $r \in \mathbb{N}_0$ be such that $g_1,...,g_m$ are log-analytic of order at most $r$. By Theorem A there is a definable cell decomposition $\mathcal{C}$ of $X_{\neq 0}$ such that $g_1|_C,...,g_m|_C$ are nicely $r$-log-analytically prepared in $x$ in a simultaneous way for every $C \in \mathcal{C}$. Fix $C \in \mathcal{C}$ and the corresponding center $\Theta$ for the $r$-logarithmic preparation of $g_1|_C,...,g_m|_C$.

{\bf Claim}\\
Let $l \in \{-1,...,e\}$ and $h \in \log(Q) \cup \{f\}$ be $l$-prepared with respect to $L$ and $Q$. Then $h|_C$ is $(l,r)$-prepared in $x$ with center $\Theta$ with respect to $P:=Q|_C$. %Per definition 0 \in P!!

{\bf Proof of the claim:} We do an induction on $l$. If $l=-1$ it is clear. So assume $l \geq 0$. Since $h|_C$ is $l$-prepared with respect to $\{g_1|_C,...,g_m|_C\}$ we have that 
$$h|_C=\mu e^\sigma \tilde{v}(\nu_1e^{\tau_1},...,\nu_ke^{\tau_k})$$
where $k \in \mathbb{N}$, the functions $\mu,\nu_1,...,\nu_k \in L|_C$ are nicely $r$-log-analytically prepared in $x$ with center $\Theta$, the functions $\sigma,\tau_1,...,\tau_k \in \log(P)$ are $(l-1)$-prepared with respect to $P$ and $L|_C$ and 
$\nu_j(t,x)e^{\tau_j(t,x)} \in [-1,1]$ for every $(t,x) \in C$ and $j \in \{1,...,k\}$. Additionally $\tilde{v}$ is a power series which converges absolutely on an open neighborhood of $[-1,1]^k$ with $\tilde{v}([-1,1]^k) \subset \mathbb{R}_{>0}$.
By the inductive hypothesis we see that $\sigma$ and $\tau_1,...,\tau_k$ are $(l-1,e)$-prepared in $x$ with center $\Theta$ with respect to $P$.

With composition of power series we see that there is a $C$-nice $a:\pi(C) \to \mathbb{R}$, a nice $r$-logarithmic scale $\mathcal{Y}$ with center $\Theta$, $s \in \mathbb{N}$, $q,p_1,...,p_s \in \mathbb{Q}^{r+1}$, a tuple $b=(b_1,...,b_s)$ of $C$-nice functions on $\pi(C)$, a tuple $\exp(d):=(\exp(d_1),...,\exp(d_s))$ of definable functions on $C$ with $d_j \in \{\tau_1,...,\tau_k\} \cup \{0\}$ for every $j \in \{1,...,s\}$ such that
$$h(t,x)=a(t)\vert{\mathcal{Y}(t,x)}\vert^{\otimes q} e^{\sigma(t,x)} v(\phi_1(t,x),...,\phi_s(t,x))$$
for every $(t,x) \in C$ where for $j \in \{1,...,s\}$
$$\phi_j:C \to [-1,1], (t,x) \mapsto b_j(t)\vert{\mathcal{Y}(t,x)}\vert^{\otimes p_j}\exp(d_j(t,x)),$$
$v$ is a power series which converges on an open neighborhood of $[-1,1]^s$ and $v([-1,1]^s) \subset \mathbb{R}_{>0}$. So we see that $h$ is $(l,r)$-prepared in $x$ with respect to $P$. (By (1) in Remark \ref{3.6} we have $e^0 \in Q$ since $l \geq 0$.)
\hfill$\blacksquare_{\textnormal{Claim}}$

Because $f$ is $e$-prepared with respect to $Q$ we see by the claim that $f|_C$ is $(e,r)$-prepared in $x$ with center $\Theta$ with respect to $P$. With the claim we obtain that $g|_C$ is $(l,r)$-prepared in $x$ with center $\Theta$ with respect to $P$ for every $g \in \log(Q)$ which is $l$-prepared with respect to $Q$. So with $(*_e)$ we see that $(+_e)$ is satisfied.
\hfill$\blacksquare$

Theorem $C$ allows us to study specific classes of definable functions and discuss various analytic properties of them. One example are the so-called \textit{restricted log-exp-analytic functions}. These are definable functions which are compositions of log-analytic functions and exponentials of locally bounded functions (see Definition 2.5 in \cite{10.1215/00192082-10234104}). For a restricted log-exp-analytic function $f:X \to \mathbb{R}$ we find a pair $(e,r) \in (\mathbb{N}_0 \cup \{-1\}) \times \mathbb{N}_0$ and a decomposition $\mathcal{C}$ of $X$ into finitely many definable cells such that for every $C \in \mathcal{C}$ there is a set $P$ of locally bounded functions on $C$ with respect to $X$ such that $f|_C$ is $(e,r)$-prepared with respect to $P$. Consequences of this observation are the following. On the one hand differentiability properties like strong quasianalticity or a parametric version of Tamm's theorem could be established for restricted log-exp-analytic functions (see \cite{10.1215/00192082-10234104}) and on the other hand it could be shown that a real analytic restricted log-exp-analytic function has a holomorphic extension which is again restricted log-exp-analytic (see \cite{OprisComplex}). 
(These are generalizations of the corresponding results from \cite{Tamm} and \cite{KaiserComplex} in the globally subanalytic setting.) So restricted log-exp-analytic functions share their properties with globally subanalytic ones from the point of analysis, but not from the point of o-minimality.

I am hopeful that these results and proofs may serve as stepping stones towards a better understanding of the o-minimal structure $\mathbb{R}_{\text{an,exp}}$ by identifying more different interesting classes of definable functions and proving analytic properties of them.

\bibliographystyle{jloganal}
\bibliography{references}

% For the article "O-minimal preparation theorems there is no doi on researchgate/similar pages
%For the article "A preparation theorem for Weiherstrass systems2 the doi in research gate is not correct! Therefore I wrote only a url.
%The article "Global complexification of restricted log-exp-analytic functions is still a preprint on arxiv.

\vs{1cm}
Andre Opris\\
University of Passau\\
Faculty of Computer Science and Mathematics\\
andre.opris@uni-passau.de\\
D-94030 Germany
\end{document}